\documentclass[sn-mathphys]{sn-jnl} %Math and Physical Sciences Reference Style
\usepackage{natbib}
\usepackage{mathtools}
\jyear{2021}%

\usepackage{lineno}

\theoremstyle{thmstyleone}%
\newtheorem*{theoremunnumbered}{Theorem}%  
\newtheorem*{remark}{Remark}

\theoremstyle{thmstyletwo}%

\theoremstyle{thmstylethree}%

\begin{document}

%\linenumbers

\title[Calculus ++ ]{\begin{center}Calculus ++\end{center} Generalized Differential and Integral Calculus and Heisenberg’s Uncertainty Principle}


\author*[1]{\fnm{Fernando} \sur{Marques de Almeida Nogueira}}\email{fernando@engenharia.ufjf.br}



\affil*[1]{\orgdiv{Department of Industrial and Mechanical Engineering}, \orgname{Juiz de Fora Federal University}, \orgaddress{\street{University Campus}, \city{Juiz de Fora}, \postcode{36036-900}, \state{MG}, \country{Brazil}}}

%\affil[2]{\orgdiv{Department}, \orgname{Organization}, \orgaddress{\street{Street}, \city{City}, \postcode{10587}, \state{State}, \country{Country}}}

%\affil[3]{\orgdiv{Department}, \orgname{Organization}, \orgaddress{\street{Street}, \city{City}, \postcode{610101}, \state{State}, \country{Country}}}


\abstract{This paper presents a generalization for Differential and Integral Calculus. 
Just as the derivative is the instantaneous angular coefficient of the tangent line to a function, the generalized derivative is the instantaneous parameter value of a reference function (derivator function) tangent to the function.
The generalized integral reverses the generalized derivative, and its calculation is presented without antiderivatives.
Generalized derivatives and integrals are presented for polynomial, exponential and trigonometric derivators and integrators functions.
As an example of the application of Generalized Calculus, the concept of instantaneous value provided by the derivative is used to precisely determine time and frequency (or position and momentum) in a function (signal or wave function), opposing Heisenberg's Uncertainty Principle.}

\keywords{Differential and Integral Calculus, Instantaneous Frequency, Heisenberg's Principle Uncertainty}



\maketitle

\section{Introduction}\label{sec:INTRODUCAO}

%The Introduction section, of referenced text \cite{bib1} expands on the background of the work (some overlap with the Abstract is acceptable). The introduction should not include subheadings.

%Springer Nature does not impose a strict layout as standard however authors are advised to check the individual requirements for the journal they are planning to submit to as there may be journal-level preferences. When preparing your text please also be aware that some stylistic choices are not supported in full text XML (publication version), including coloured font. These will not be replicated in the typeset article if it is accepted.

Differential and Integral Calculus has become one of the main mathematical tools that made possible discoveries and advances in several areas such as Physics, Chemistry, Economics, Computer Science, Engineering, and even Biology and Medicine. Moreover, in Mathematics itself, Differential and Integral Calculus is used in other areas, such as Linear Algebra, Analytical Geometry, Probability, and Optimization, among others.

Differential and Integral Calculus was developed by Isaac Newton \cite{Newton1736} (1643-1727) and Gottfried Wilhelm Leibniz \cite{Leibniz1684} (1646-1716), independently of each other, in the 17th century and basically established three operations that are applicable to any function: calculus of limits, derivatives, and integrals. 

The derivative concerns the instantaneous rate of change of a function. On the other hand, the integral concerns the area under the curve described by a function. Both the derivative and the integral are based on the calculus of infinitesimals through the concept of limit, and the Fundamental Theorem of Calculus formalizes the inverse operations relationship between Differential and Integral Calculus.

The derivative is an operation that is performed on any function $f(x)$\footnote{$f(x)$ is formally defined in the remaining sections.}, resulting in another function $f'(x)$ that represents the slope of the tangent line to $f(x)$ for each x.
Differential Calculus uses the line as the ``reference function" and its slope as the result of the derivative. 
\begin{remark}
\textbf{Why use only the line as the reference function and its slope as the result of the derivative?}
\end{remark}
This paper presents the derivative performed for other reference functions different from the line and other parameters different from the slope of the line, thus generalizing the Differential Calculus.

Since the derivative and the integral are inverse operations, the same generalization concept employed for Differential Calculus is applied to Integral Calculus.

\subsection{The Derivative, its Generalization and the Antiderivative}\label{subsec:DERIVATIVE}

\begin{figure}[h]%
\centering
\includegraphics[width=0.9\textwidth]{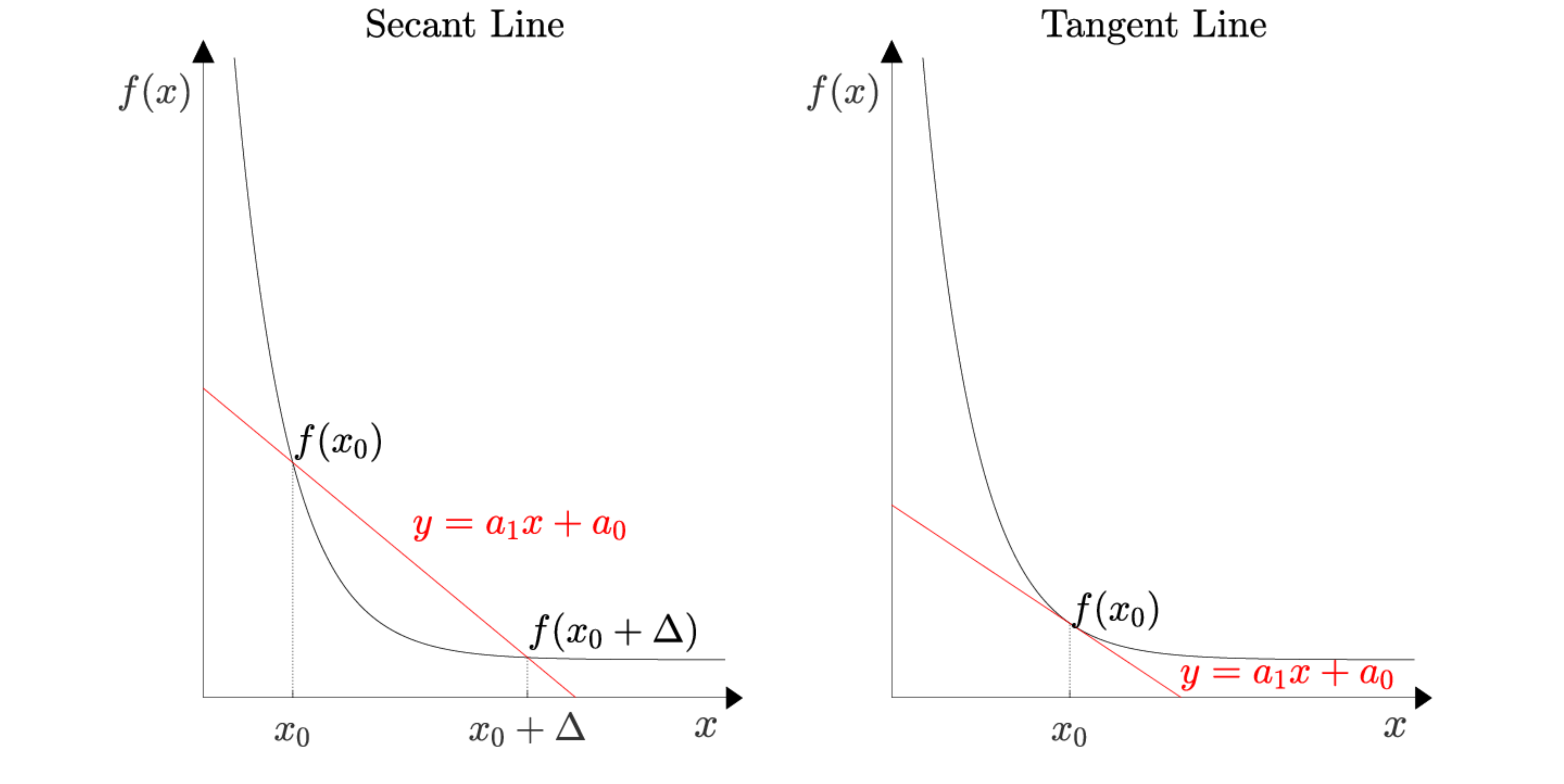}
\caption{Secant line (red) passing through the points $(x_0,f(x_0))$ and $(x_0+\Delta,f(x_0+\Delta))$ (left figure) and the tangent line at point $(x_0,f(x_0))$ (right figure) to the function (black).}\label{secante_line}
\end{figure}

The derivative of a function can be understood as a linear interpolation process. Let $\mathbb{I}$ a non-empty open interval, $f:\mathbb{I} \to \mathbb{R}$ a function, $y=f(x)$,  $\mathbb{I} \subseteq \mathbb{R}$, $x_0 \in \mathbb{I}$ and $\Delta \in \mathbb{R}$, as illustrated in figure \ref{secante_line}.

Two points determine a line: from the points $(x_0,f(x_0))$ and $(x_0+\Delta,f(x_0+\Delta))$ it is possible to calculate the angular ($a_1$) and linear ($a_0$) coefficient of the linear equation $y = a_1x + a_0$ secant to the graph of the function $f(x)$. This calculation is obtained by solving the following linear system:

\begin{equation} \label{sistema_linear1}
S:\left\{
\begin{array}{ll}
f(x_0)=a_1x_0+a_0\\
f(x_0+\Delta)=a_1(x_0+\Delta)+a_0
\end{array}
\right.
\end{equation}

The resolution of  \ref{sistema_linear1} is:
\begin{equation} \label{a_solucao_sistema_linear1}
a_1= \frac{f(x_0+\Delta)-f(x_0)}{\Delta}
\end{equation}

\begin{equation} \label{b_solucao_sistema_linear1}
a_0=\frac{f(x_0)(x_0+\Delta)-f(x_0+\Delta)x_0}{\Delta}
\end{equation}

In differential calculus, the angular coefficient ($a_1$)  in  \eqref{a_solucao_sistema_linear1} is known as Newton's Difference Quotient.

For small $\Delta$ values, the linear equation $y = a_1x + a_0$ will be practically tangential to the graph of the function $f(x)$ near the point $x_0$, and in the limit $\Delta \to 0$, this line will be tangential to the graph of $f(x)$ at point $x_0$.
Applying limit of $\Delta \to  0$ in \eqref{sistema_linear1} is:

\begin{equation} \label{Limite_sistema_linear1}
S:\left\{
\begin{array}{ll}
f(x_0)=a_1x_0+a_0\\
f(x_0+\Delta \to 0)=a_1(x_0+\Delta \to 0)+a_0
\end{array}
\right.
\end{equation}

And, its resolution is:
\begin{equation} \label{a_Limitesolucao_sistema_linear1}
a_1\mid_{x_0}= \displaystyle \lim_{\Delta \to 0} \frac{f(x_0+\Delta)-f(x_0)}{\Delta}
\end{equation}

\begin{equation} \label{b_Limitesolucao_sistema_linear1}
a_0\mid_{x_0}= \displaystyle \lim_{\Delta \to 0} \frac{f(x_0)(x_0+\Delta)-f(x_0+\Delta)x_0}{\Delta}
\end{equation}

In \eqref{a_Limitesolucao_sistema_linear1}, the value of $a_1\mid_{x_0}$ is the value of the derivative of $f(x)$ at the point $x_0$. The value of $a_0\mid_{x_0}$ in \eqref{b_Limitesolucao_sistema_linear1} \textbf{is not used} in traditional differential and integral calculus. Generalizing for any point $x$ in the domain, $a^{ins}_1: \mathbb{I} \to \mathbb{R}$ the function $a_1$ instantaneous, the derivative of $f(x)$ is:

\begin{equation} \label{derivadaGeral}
\frac{df(x)}{dx}=a^{ins}_1(x)= \lim_{\Delta \to 0} \frac{f(x+\Delta)-f(x)}{\Delta}
\end{equation}

\begin{remark}
Therefore, the derivative of a function is the application of the limit of $\Delta \to  0$ to Newton's Difference Quotient.
\end{remark}
This process uses a linear procedure to determine the \textbf{slope} (angular coefficient) of the linear equation of the tangent line to any function at a given point. This linear procedure is simply a linear interpolation or a regression to the linear equation for two infinitesimally close points belonging to $f(x)$.
However, similarly to \eqref{derivadaGeral}, from \eqref{b_Limitesolucao_sistema_linear1}, $a^{ins}_0: \mathbb{I} \to \mathbb{R}$ the function $a_0$ instantaneous,  one can write:

%$\mathfrak{D}\{a_0\} = a_1x+a_0$
\begin{equation} \label{derivadaGeralParametro b}
\begin{array}{ll}
\displaystyle\mathfrak{D}\{\}= a_1x+a_0\\
\displaystyle\mathfrak{D}\{a_0\}\frac{df(x)}{dx}=a^{ins}_0(x)= \lim_{\Delta \to 0} \frac{f(x)(x+\Delta)-f(x+\Delta)x}{\Delta}
\end{array}
\end{equation}

where,

$\mathfrak{D}\{\}$ is the tangent function to $f(x)$ in which the derivative is defined (the line equation in this case, as in the ``classical derivative", but it could be any other function);

$\mathfrak{D}\{a_0\}\frac{df(x)}{dx}$ indicates under which parameter of the $\mathfrak{D}\{\}$ the derivative $\frac{df(x)}{dx}$ is defined.

In this context, as \eqref{derivadaGeralParametro b}, the generalized notation for the ``classical derivative" \eqref{derivadaGeral} is:
\begin{equation} \label{derivadaGeralParametro a}
\begin{array}{ll}
\displaystyle\mathfrak{D}\{\}= a_1x+a_0\\
\displaystyle\mathfrak{D}\{a_1\}\frac{df(x)}{dx}=a^{ins}_1(x)= \lim_{\Delta \to 0} \frac{f(x+\Delta)-f(x)}{\Delta}
\end{array}
\end{equation}

The reasoning used in \eqref{sistema_linear1} to \eqref{derivadaGeralParametro a} can be generalized to other functions (and not just the linear equation) and their respective parameters.
This concept can also be applied to the integral of a function. For example, the notation for the inverse operation of \eqref{derivadaGeralParametro a} is:
\begin{equation} \label{intregralGeralParametro a}
\begin{array}{ll}
\displaystyle\mathfrak{I}\{\}= a_1x+a_0\\ 
\displaystyle F(x)=\mathfrak{I}\{a_1\}\int f(x)dx
%F(x)= \(\int f(x)\,dx[ax+\underline{b}]\)
\end{array}
\end{equation}

where,

$\mathfrak{I}\{\}$ is the tangent function to $F:\mathbb{I} \to \mathbb{R}$ which the integral is defined;

$\mathfrak{I}\{a_1\}\int f(x)dx$ indicates under which parameter of the $\mathfrak{I}\{\}$ the integral of $f(x)$ is defined.

As a result, the following concepts can be defined:

\begin{itemize}
    \item \textbf{Derivator Function} is the function $\mathfrak{D}\{\}$ that is used in the interpolation process (the classic derivative uses the linear equation as Derivator Function).
    \item \textbf{Derivator Parameter} is the parameter of interest $p_k$ of the  $\mathfrak{D}\{\}$, represented by $\mathfrak{D}\{p_k\}$, where $p_k \in \mathfrak{P}$, $\mathfrak{P}=\{p_0,p_1,p_2,...,p_{N-1}\}$ is the set of $N\in \mathbb{N}$ parameters of the $\mathfrak{D}\{\}$, $k\in \mathbb{N}$ and $k<N$ . 
    (the classic derivative uses the angular coefficient $a_1$ as Derivator Parameter).
    
    \item \textbf{Integrator Function} is the function $\mathfrak{I}\{\}$ that is used in the process of obtaining the primitive function (the classic integral uses the linear equation as Integrator Function).
    \item \textbf{Integrator Parameter} is the parameter of interest $p_k$ of the  $\mathfrak{I}\{\}$, represented by $\mathfrak{I}\{p_k\}$, where $p_k \in \mathfrak{Q}$, $\mathfrak{Q}=\{p_0,p_1,p_2,...,p_{N-1}\}$ is the set of $N\in \mathbb{N}$ parameters of the $\mathfrak{I}\{\}$, $k\in \mathbb{N}$ and $k<N$ . 
    (the classic integral uses the angular coefficient $a_1$ as Integrator Parameter).
\end{itemize}

%Similarly to the Derivator Function and the Derivator parameter, in \ref{intregralGeralParametro b}, [\textit{$a_1x$}+\textit{\underline{$a_0$}}] is the \textbf{Integrator Function} and ($a_0$) is the \textbf{Integrator Parameter}.

\begin{figure}[h] 
\centering
\includegraphics[width=1\textwidth]{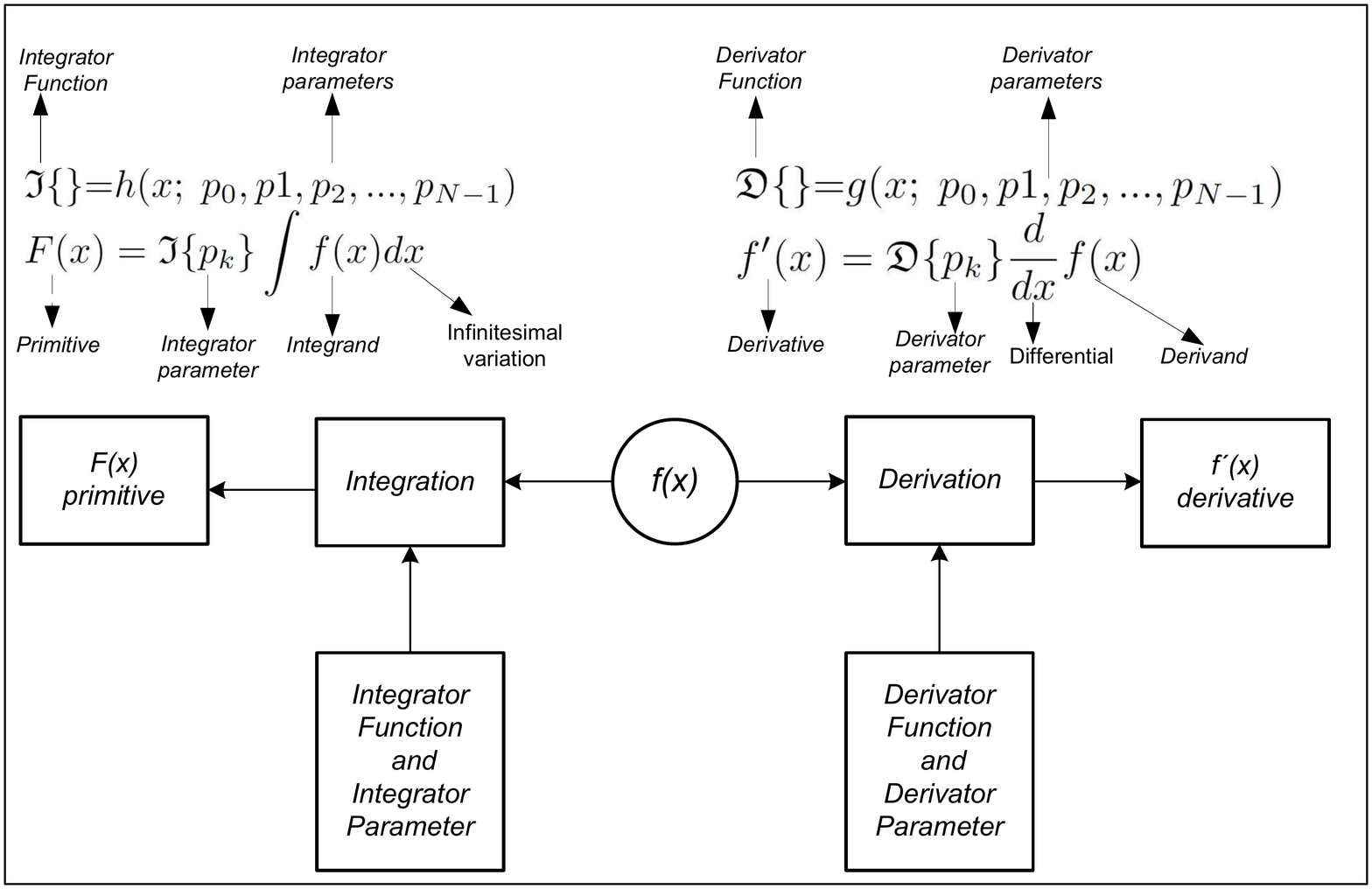}
\caption{Functions and operations involved in Differential and Integral Calculus.}\label{figureCalculus++}
\end{figure}

Figure \ref{figureCalculus++} illustrates the names of the functions and operations involved in Generalized Differential and Integral Calculus.

\section{Background}\label{sec:BACKGROUND}

Different forms of the derivative have already been established. These forms use concepts different from the foundation employed for generalizing Differential and Integral Calculus presented in this article.
\begin{itemize}
\item Symmetric Derivative \cite{Aull1967}
\end{itemize}
A simple variant form of the ``classical derivative" is the Symmetric Derivative, which uses Newton's Difference Quotient in a symmetrical form. The Symmetric Derivative $f'_S$ is defined as:

\begin{equation} \label{derivadaSimetrica}
f'_S=\lim_{\Delta \to 0} \frac{f(x+\Delta)-f(x-\Delta)}{2\Delta}
\end{equation}

Although the Symmetric Derivative uses a different form for Newton's Difference Quotient, the derivative function can still be understood as the slope of the tangent line to the function $f(x), \forall{x}$. In this form, the Symmetric Derivative is a different manner of defining the ``classical derivative".

\begin{itemize}
\item Fréchet Derivative \cite{Coleman2012}
\end{itemize}
Given $\mathbb{V}$ and $\mathbb{U}$ normed vectorial spaces, $\mathbb{W} \subseteq \mathbb{V}$, $f:\mathbb{W} \to \mathbb{U}$ a function Fréchet differentiable at $x \in \mathbb{V}$. If there is a bounded and linear operator $A:\mathbb{V} \to \mathbb{U}$ such that:

\begin{equation} \label{derivadaFrechet}
\lim_{\Delta \to 0} \frac{\parallel f(x+\Delta)-f(x)-A\Delta\parallel_\mathbb{U}}{\parallel\Delta\parallel_\mathbb{V}} = 0
\end{equation}

Then, $A$ is the derivative of $f$ at $x$.
The Fréchet Derivative is used on a vector-valued function of multiple real variables and to define the Functional Derivative, generalizing the derivative of a real-valued function of a single real variable.

\begin{itemize}
\item Functional Derivative \cite{Coleman2012}
\end{itemize}

Another form of a derivative is the Functional Derivative. 
Given $\mathbb{V}$ a vectorial (function) space, $\mathbb{K}$ a field and $F$ a functional,  $F:\mathbb{V} \to \mathbb{K}$, $f\in \mathbb{V}$, $\zeta$ an arbitrary function, the Functional Derivative of $F$ at $f$, $\frac{d(F)}{d(f)}$ is:

\begin{equation} \label{derivadaFuncional}
\int{}{}{\frac{d(F)}{d(f)}}(x)\zeta(x)dx=\lim_{\Delta \to 0} \frac{F(f+\Delta\zeta)-F(f)}{\Delta}
\end{equation}

In this case, the concept of the derivative is applied to a functional and not to a function. In this paper, the concept of the derivative is generalized to functions.

\begin{itemize}
\item Fractional Derivative \cite{Golmankhaneh2022}
\end{itemize}
The derivative can be repeated $n$ times over a function, resulting in the derivative's order. Thus, the order of the derivative is clearly a natural number ($n \in \mathbb{N})$.
The fractional derivative generalizes the concept of derivative order so that the $\alpha$ order of the fractional derivative is $\alpha \in \mathbb{R}$ or even $\alpha \in \mathbb{C}$.
It is then possible, under this generalization, to calculate the derivative of $f(x)$ of order $alpha = 2.5$ or $alpha = -1$ (integral of $f(x)$), for example. In this paper, the derivative is generalized to functions and not to the order of the derivative.

\begin{itemize} 
\item q-Derivative \cite{Chaundy1962}
\end{itemize}

The q-Derivative of a function $f(x)$ is a q-analog of the ``classic derivative". Let $q \in \mathbb{R}$, it is given by:
\begin{equation} \label{Qderivada}
{\left(\frac{d}{dx}\right)_q f(x)}= \frac{f(qx)-f(x)}{qx-x}
\end{equation}
For $q\to 1$, the q-Derivative is the ``classic derivative".

\begin{itemize}
\item Arithmetic Derivative \cite{Haukkanen2018}
\end{itemize}
Let $a.b \in \mathbb{N}$ and $p$ a prime number, the arithmetic derivative $D(a.b)$ is such that:
\begin{equation} \label{Arithmetic_derivative}
\begin{array}{ll}
D(0)=D(1)=0\\
D(p)=1, \forall \mbox{ prime }p\\
D(a.b)=D(a)b+a.D(b)
\end{array}
\end{equation}
The Arithmetic Derivative is a "number derivative", which is based on prime factorization.
The arithmetic derivative can be extended to rational numbers.

Other forms of derivatives include:
\begin{itemize}
\item Carlitz derivative \cite{Kochubei2007}
\item Covariant derivative \cite{Zuk1986}
\item Dini derivative \cite{Kannan1996}
\item Exterior derivative \cite{Tu2007}
\item Gateaux derivative \cite{Andrews2011}
\item H derivative \cite{Kac2002}
\item Hasse derivative \cite{Hoffmann2015}
\item Lie derivative \cite{Lee2003}
\item Pincherle derivative \cite{Mainardi2011}
\item Quaternionic derivative \cite{Xu2015}
\item Radon Nikodym derivative \cite{Konstantopoulos2011}
\item Semi differentiability \cite{Luc2018}
\item Subderivative \cite{Rockafellar1998}
\item Weak derivative \cite{Lamboni2022}
\end{itemize}

All of these forms use concepts different from the foundation employed for the generalization of Differential and Integral Calculus presented in this article.

\section{Polynomials Derivators Functions}\label{sec:Polynomial}

%Let $n \in \mathbb{N}$, $a_i \in \mathbb{R}$, $a_i(x):\mathbb{R} \to \mathbb{R}$, $\forall{ i} \in \mathbb{N}$ and a polynomial function is a function $P(x)$ that can be written as:

Let $n \in \mathbb{N}$, $a_i \in \mathbb{R}$, $\forall{ i} \in \mathbb{N}$, $P:\mathbb{R} \to \mathbb{R}$ is a Polynomial Function if:

\begin{equation} \label{polinomio}
P(x)= a_nx^n+a_{n-1}x^{(n-1)}+....+a_1x^1+a_0x^0=\sum_{i=0}^{n} a_ix^i  \end{equation}

\begin{figure}[h] \label{figureDerivadaParabolica}
\centering
\includegraphics[width=0.9\textwidth]{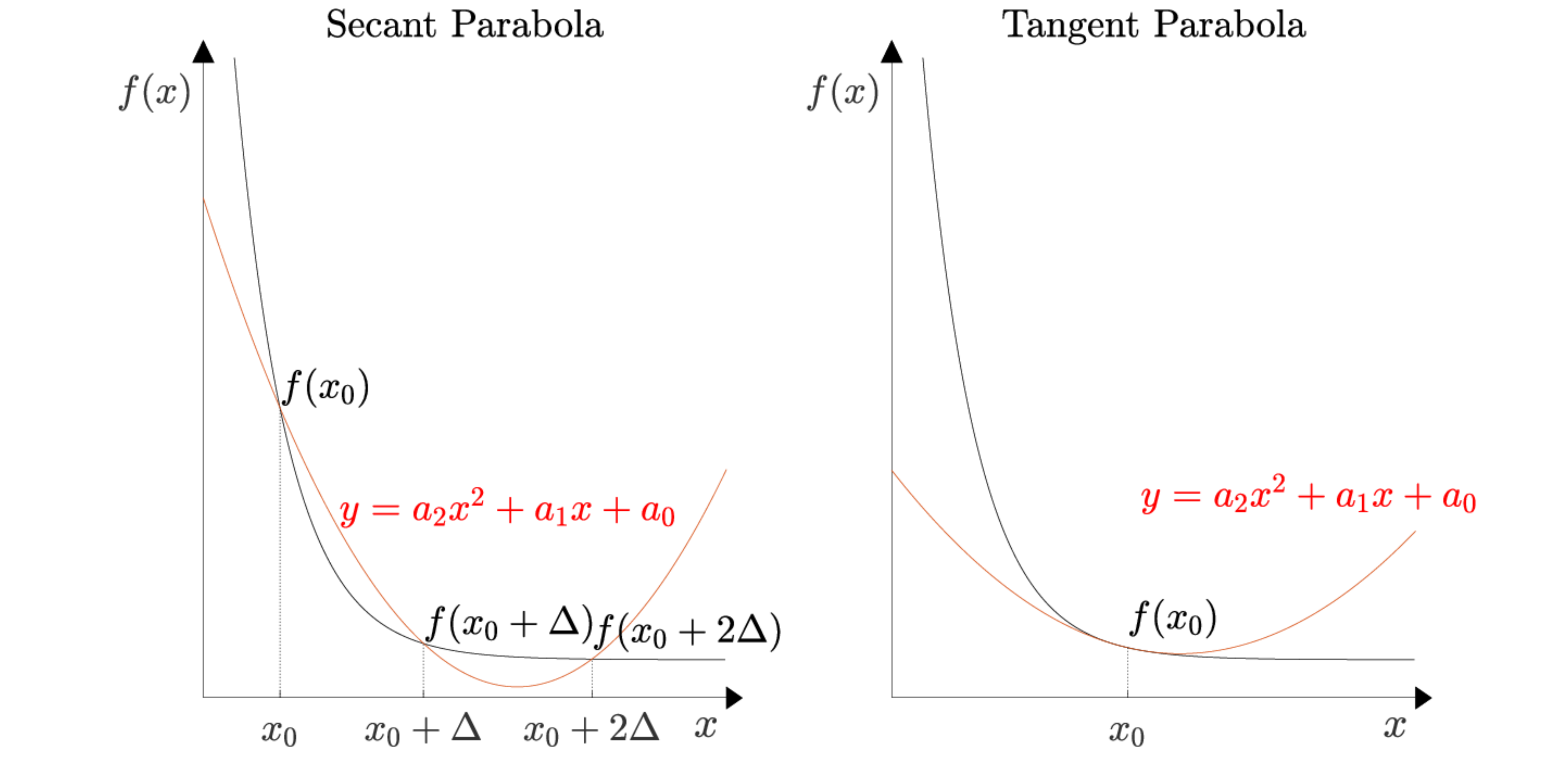}
\caption{Secant parabola (red) passing through the points $x_0$, $x_0+\Delta$ and $x_0+2\Delta.$ (left figure) and tangent parabola at point $x_0$ (right figure) to the function (black).  }\label{fig3}
\end{figure}

%In \ref{polinomio}, the value of n (with a_n 0,obviously) defines %the degree of the polynomial.
In \eqref{polinomio}, the value of $n$ defines the degree of the polynomial.
For $n=1$, the polynomial is a linear equation, and two points are needed to define its parameters (as in \eqref{sistema_linear1}). For $n=2$ and $n=3$, the polynomial is a quadratic (parabola) and cubic function, and 3 and 4 points are needed to define their parameters, respectively. For other degrees of the polynomial, the reasoning is analogous; therefore, $n+1$ points are necessary to define the parameters of a polynomial of degree n.
Using a polynomial of degree 2 ($n=2$) as derivator function, the derivative becomes an interpolation process to the quadratic function for three infinitesimally close points belonging to $f(x)$, resulting in the \textbf{Parabolic Derivative}, as shown in figure \ref{figureDerivadaParabolica}.

The system is:

\begin{equation} \label{sistema_parabola}
S:\left\{
\begin{array}{ll}
f(x_0)=a_2x_0^2+a_1x_0+a_0\\
f(x_0+\Delta)=a_2(x_0+\Delta)^2+a_1(x_0+\Delta)+a_0\\
f(x_0+2\Delta)=a_2(x_0+2\Delta)^2+a_1(x_0+2\Delta)+a_0
\end{array}
\right.
\end{equation}

Generalizing for any point $x$ in the domain, $a^{ins}_0: \mathbb{I} \to \mathbb{R}$ the function $a_0$ instantaneous, $a^{ins}_1: \mathbb{I} \to \mathbb{R}$ the function $a_1$ instantaneous, $a^{ins}_2: \mathbb{I} \to \mathbb{R}$ the function $a_2$ instantaneous, and applying limit to $\Delta \to 0$, the resolution of \eqref{sistema_parabola} is:

\begin{equation} \label{derivadoraParabolica}
\mathfrak{D}\{\}=a_2x^2+a_1x+a_0
\end{equation}

\begin{eqnarray} \label{a2_solucao_sistema_parabola}
\mathfrak{D}\{a_2\}\frac{df(x)}{dx}=a^{ins}_2(x)=\lim_{\Delta \to 0} \frac{1}{2} \frac{f(x)-2f(x+\Delta)+f(x+2\Delta)}{\Delta^2}
\end{eqnarray}

\begin{equation} \label{a1_solucao_sistema_parabola} 
\mathfrak{D}\{a_1\}\frac{df(x)}{dx}=a^{ins}_1(x)=      
\lim_{\Delta \to 0} -\frac{1}{2} \frac{K_0f(x)-K_1f(x+\Delta)+K_2f(x+2\Delta)}{\Delta^2}\\
\end{equation}
where,
\begin{align}
K_0 &= 2x+3\Delta \nonumber \\
K_1 &= 4x+4\Delta \nonumber \\
K_2 &= 2x+\Delta \nonumber 
\end{align}

%\begin{eqnarray} \label{a1_solucao_sistema_parabola}
%\mathfrak{D}\{a_1\}\frac{df(x)}{dx}=...\\
%\lim_{\Delta \to 0} -\frac{1}{2} \frac{(2x+3\Delta)f(x)-(4x+4\Delta)f(x+\Delta)+(2x+\Delta)f(x+2\Delta)}{\Delta^2}
%\end{eqnarray}

\begin{equation} \label{a0_solucao_sistema_parabola}
\mathfrak{D}\{a_0\}\frac{df(x)}{dx}=a^{ins}_0(x)=
\lim_{\Delta \to 0} \frac{1}{2} \frac{K_0f(x)+K_1f(x+\Delta)+K_2f(x+2\Delta)}{\Delta^2}\\
\end{equation}
where,
\begin{align}
K_0 &= 2\Delta^2+x^2+3x\Delta \nonumber \\
K_1 &= -2x^2-4x\Delta \nonumber \\
K_2 &= x^2+x\Delta  \nonumber
\end{align}

%\begin{eqnarray} 
%\label{a0_solucao_sistema_parabola}
%\mathfrak{D}\{a_0\}\frac{df(x)}{dx}=...\\
%\lim_{\Delta \to 0} \frac{1}{2} \frac{(2\Delta^2+x^2+3x\Delta)f(x)+(-2x^2-%4x\Delta)f(x+\Delta)+(x^2+x\Delta)f(x+2\Delta)}{\Delta^2}
%\end{eqnarray}

% \begin{eqnarray}
% &\mathfrak{D}\{\}:a_2x^2+a_1x+a_0 \label{derivadoraParabolica} \\
% \mathfrak{D}\{a_2\}\frac{df(x)}{dx}&=\lim_{\Delta \to 0} \frac{1}{2} \frac{f(x)-2f(x+\Delta)+f(x+2\Delta)}{\Delta^2} \\ \label{a2_solucao_sistema_parabola}
% \mathfrak{D}\{a_1\}\frac{df(x)}{dx}&=\lim_{\Delta \to 0} -\frac{1}{2} \frac{(2x+3\Delta)f(x)-(4x+4\Delta)f(x+\Delta)+(2x+\Delta)f(x+2\Delta)}{\Delta^2} \label{a1_solucao_sistema_parabola}\\ 
% \mathfrak{D}\{a_0\}\frac{df(x)}{dx}&=\lim_{\Delta \to 0} \frac{1}{2} \frac{(2\Delta^2+x^2+3x\Delta)f(x)+(-2x^2-4x\Delta)f(x+\Delta)+(x^2+x\Delta)f(x+2\Delta)}{\Delta^2} \label{a0_solucao_sistema_parabola} \\
% \end{eqnarray}

In this form, \eqref{a2_solucao_sistema_parabola}, \eqref{a1_solucao_sistema_parabola} and \eqref{a0_solucao_sistema_parabola} define the Parabolic Derivative to the $a_2$, $a_1$ and $a_0$ parameters, respectively. For polynomials of other degrees, the procedure is similar to that performed in \eqref{sistema_parabola}. 
For example, the generalized polynomial derivative (polynomial derivator function) of $f(x)=cx^m$, $c,m \in \mathbb{R}$ for the highest degree parameter $n \in \mathbb{N^*}$ of the polynomial derivator function is:

\begin{equation} \label{derivadora_generalizada_polinomial}
\mathfrak{D}\{\} = {a_n}x^n+a_{n-1}x^{(n-1)}+...+a_1x^1+a_0x^0
\end{equation}

\begin{eqnarray} \label{derivada_generalizada_polinomial}
\mathfrak{D}\{{a_n}\}
\frac{d(cx^m)}{dx}=c\left(\prod_{i=0}^{n-1} \frac{m-i}{i+1}\right)x^{(m-n)}=\; & & \nonumber \\
\frac{c}{n!}\left(\prod_{i=0}^{n-1} {(m-i)}\right)x^{(m-n)}\;\;\;\;\;
\end{eqnarray}

The Antiderivative of \eqref{derivada_generalizada_polinomial} is:

\begin{equation} \label{integradora_generalizada_polinomial}
\mathfrak{I}\{\}=\mathfrak{D}\{\} = {a_n}x^n+a_{n-1}x^{(n-1)}+...+a_1x^1+a_0x^0
\end{equation}
\begin{equation} \label{antiderivada_generalizada_polinomial}
\begin{array}{ll}
\mathfrak{I}\{{a_n}\}\int{cx^m}{dx}=\left(\frac{cx^{(m+n)}}{(m+n)\prod_{i=1}^{n-1}{\frac{m+i}{i+1}}}\right),\forall{(m+n)\neq{0} \wedge n \geq 2}\\
\mathfrak{I}\{{a_1}\}\int{cx^m}{dx}=\left(\frac{cx^{(m+1)}}{(m+1)}\right),m\neq{-1}
\end{array}
\end{equation}

For the other functions $f(x)$ and/or other derivator parameters, the reasoning is analogous to \eqref{sistema_parabola}, \eqref{derivada_generalizada_polinomial} and \eqref{antiderivada_generalizada_polinomial}.

\subsection{Vanishing Terms and Primitives}\label{subsec:Vanishing}

The ``classic derivative" uses the linear equation as the derivator function and the angular coefficient as the derivator parameter (as in \eqref{derivadaGeralParametro a}). However the derivative in this form does not model the linear coefficient of the derivator function, and therefore this term, if it exists in the derivand $f(x)$, ``vanishes" for the differential operator, not influencing the derivative.
Using the linear equation as the derivator function and the linear coefficient as the derivator parameter as in \eqref{derivadaGeralParametro b}, the derivative in this form does not model the first-degree term (angular coefficient) of the derivator function. Therefore, this term "vanishes" for the differential operator, not influencing the derivative.
The following example is suitable for showing this case. Considering:
\begin{equation} \label{exemplo função}
f(x)=x^2+2x+3
\end{equation}
for $\mathfrak{D}\{\}=a_1x+a_0$, their derivatives are:

\begin{equation} \label{derivada_exemplo_funçãoParametro a}
\mathfrak{D}\{a_1\}\frac{df(x)}{dx}=2x+2
\end{equation}
and

\begin{equation} \label{derivada_exemplo_funçãoParametro b}
\mathfrak{D}\{a_0\}\frac{df(x)}{dx}=-x^2+3
\end{equation}

The term $+3$ and $2x$ in \eqref{exemplo função} vanishes in \eqref{derivada_exemplo_funçãoParametro a} and \eqref{derivada_exemplo_funçãoParametro b}, respectively. For $\mathfrak{I}\{\}=a_1x+a_0$, the antiderivatives, respectively, for \eqref{derivada_exemplo_funçãoParametro a} and \eqref{derivada_exemplo_funçãoParametro b} are:

\begin{equation} \label{antiderivada_exemplo_funçãoParametro a}
%\(\int (2x+2)\,dx[\underline{a_1}x+a_0}]\)=x^2+2x
\mathfrak{I}\{a_1\}\int (2x+2)dx=x^2+2x
\end{equation}
and
\begin{equation} \label{antiderivada_exemplo_funçãoParametro b}
\mathfrak{I}\{a_0\}\int (-x^2+3)dx=x^2+3
%\(\int (-x^2+3)\,dx[a_1x+\underline{a_0}]\)=x^2+3
\end{equation}

Since \eqref{derivada_exemplo_funçãoParametro a} and \eqref{derivada_exemplo_funçãoParametro b} do not model the terms $a_0$ and $a_1x$, respectively, the antiderivatives \eqref{antiderivada_exemplo_funçãoParametro a} and \eqref{antiderivada_exemplo_funçãoParametro b} do not return in \eqref{exemplo função} and must be added by the following terms ($k_0$ and $k_1x$, with $k_0$ and $k_1$ constants):

\begin{equation} \label{antiderivada_exemplo_funçãoParametro a + k}
 %\(\int (2x+2)\,dx[\underline{a_1}x+a_0}]\)=x^2+2x+k
\mathfrak{I}\{a_1\}\int (2x+2)dx=x^2+2x+k_0
\end{equation}
and
\begin{equation} \label{antiderivada_exemplo_funçãoParametro b + k}
%\(\int (-x^2+3)\,dx[a_1x+\underline{a_0}]\)=x^2+kx+3
\mathfrak{I}\{a_0\}\int (-x^2+3)dx=x^2+k_1x+3
\end{equation}

The addition of the terms $k_0$ and $k_1x$ in \eqref{antiderivada_exemplo_funçãoParametro a + k} and \eqref{antiderivada_exemplo_funçãoParametro b + k} is necessary because, \textbf{independently of $k_0$ and $k_1x$}, their derivatives are the same:

\begin{equation} \label{derivada_exemplo_funçãoParametro a com k}
\mathfrak{D}\{a_1\}\frac{d(x^2+2x+k_0)}{dx}=2x+2
\end{equation}
and

\begin{equation} \label{derivada_exemplo_funçãoParametro b com k}
\mathfrak{D}\{a_0\}\frac{d(x^2+k_1x+3)}{dx}=-x^2+3
\end{equation}

\subsection{Integrals without Antiderivatives}\label{subsec:Integrals}
The Fundamental Theorem of Calculus (FTC) \cite{Strang1991} establishes the relationship between differential calculus and integral calculus, as \textbf{inverse operations} (with reservations). The FTC is divided into two parts. Part 1 shows that the derivative of the integral of $f(x)$ is equal to $f(x)$: this is perfect! Part 2 reverses the order, that is, the integral of the derivative of $f(x)$ is equal to $f(x)$, however, plus a constant $k$, that is, $f(x)+k$: this is perfect too, but the exact return to function $f(x)$ does not occur when the derivative is performed first and then the integral. Thus, in formal terms, the FTC states that the operations of derivation and integration are inverse, apart from a constant value.

\begin{equation} \label{TFC}
\frac{d\int f(x)dx}{dx}=f(x) \neq \int \frac{df(x)}{dx}dx=f(x)+k
%\int (-x^2+3)dx[a_1x+\underline{a_0}])=x^2+kx+3
\end{equation}

This problem occurs simply because the ``classic derivative" only gives the instantaneous rate of change of a function for its domain. This rate, as seen in \eqref{a_Limitesolucao_sistema_linear1}, is the angular coefficient for the linear equation when used as derivator function. Obviously, the linear equation cannot be defined by its angular coefficient $a_1$ alone. The linear coefficient $a_0$ also needs to be defined for the linear equation to be complete.

The derivation process is carried out by applying the concept of limit to Newton's quotient. On the other hand, the integration process does not have a specific form, and this is obtained, in practice, through the calculation of antiderivatives. Nonetheless, a function can be defined by applying the integrator function $\mathfrak{I}\{\}$ to the generalized derivatives for a given derivator function $\mathfrak{D}\{\}$ for \textbf{all} their respective parameters, with $\mathfrak{I}\{\}=\mathfrak{D}\{\}$.

For the integrator function $\mathfrak{I}\{\}:a_1x+a_0$ (linear equation) and $\mathfrak{I}\{\}=\mathfrak{D}\{\}$, the primitive $f(x)$ is:

\begin{equation} \label{primitva_sem_integracao_linear}
f(x)=\mathfrak{D}\{a_1\}\frac{df(x)}{dx}x+\mathfrak{D}\{a_0\}\frac{df(x)}{dx}
%f(x)=\frac{df(x)}{dx[\underline{a_1}x+a_0]}x+\frac{df(x)}{dx[a_1x+\%underline{a_0}]}
\end{equation}

It is important to emphasize that $f(x)$ was obtained from its generalized derivatives without the conventionally used integration process (antiderivative) in ``classical integral calculus".

For example, from the functions \eqref{derivada_exemplo_funçãoParametro a} and \eqref{derivada_exemplo_funçãoParametro b} (derivatives of $f(x)=x^2+2x+3$), the primitive $f(x)$ is:
\begin{equation} \label{exemplo_primitva_sem_integracao_linear}
f(x)=(2x+2)x-x^2+3=x^2+2x+3
\end{equation}

The \textbf{exact} return to $f(x)$ is obtained (\eqref{exemplo_primitva_sem_integracao_linear} equals \eqref{exemplo função}). The concept involved in obtaining the function $f(x)$ is: 

\begin{theoremunnumbered} [{}]
Let $\mathbb{I}$ a non-empty open interval, $\mathbb{I} \subseteq \mathbb{C}$, $h:\mathbb{I} \to \mathbb{C}$ a function, $y=h(x;p_0,p1,p_2,...,p_{N-1})$ and $\mathfrak{P}=\{p_0,p_1,p_2,...,p_{N-1}\}$ the set of $N\in \mathbb{N}$ parameters. If S is a system that has a unique solution for $N$ points $(x_k,y_k)$, $k\in \mathbb{N}$, $k\leq{N-1}$, such as:
\begin{equation} \label{sistema_S}
S:\left\{
\begin{array}{ll}
y_0=h(x_0;p_0,p1,p_2,...,p_{N-1})\\
y_1=h(x_1;p_0,p1,p_2,...,p_{N-1})\\
y_2=h(x_2;p_0,p1,p_2,...,p_{N-1})\\
...\\
y_{N-1}=h(x_{N-1};p_0,p1,p_2,...,p_{N-1})\\
\end{array}
\right.
\end{equation}
$\mathfrak{I}\{\}=\mathfrak{D}\{\}=h(x;p_0,p1,p_2,...,p_{N-1})$, $f(x)$ differentiable on $\mathfrak{D}\{\}$, $\forall{x}\in \mathbb{I}$, then f(x) can be described by $\mathfrak{I}\{\}$ whose parameters are given by their generalized derivatives in their respective $N$ parameters, i.e.  $f(x)=h(x;\mathfrak{D}\{p_0\},\mathfrak{D}\{p_1\},\mathfrak{D}\{p_2\},...,\mathfrak{D}\{p_{N-1}\})$.
\end{theoremunnumbered}

 %\begin{equation} \label{Geral_primitva_sem_integracao}
 %f(x)=h\left(x,\frac{df(x)}{dx[g(x,\underline{p_0},p_1,...,p_k]},\frac{df(x)}{dx[g(x,p_0,\underline{p_1},...,p_k]},...,\frac{df(x)}{dx[g(x,p_0,p_1,...,\underline{p_k}]}\right)
 %\end{equation}

\section{Exponential Derivators Functions}\label{sec:Exponential}

Let $A,a,b,x \in \mathbb{R}$, $a^{ins}:\mathbb{R} \to \mathbb{C}$, $b^{ins}:\mathbb{R} \to \mathbb{C}$ and $f:\mathbb{R} \to \mathbb{R}$ an exponential function as:

\begin{equation} \label{exponencial1}
f(x)= Ae^{ax},A\geq0 \end{equation}

Making $A=e^b$, \eqref{exponencial1} is:
\begin{equation} \label{exponencial2}
f(x)= e^{b}e^{ax}=e^{ax+b} \end{equation}

The following system can be written:

\begin{equation} \label{sistema_exponencial}
S:\left\{
\begin{array}{ll}
f(x)=e^{ax+b}\\
f(x+\Delta)=e^{a(x+\Delta)+b}
\end{array}
\right.
\end{equation}

Solving the system and applying the limit of $\Delta \to  0$ in \eqref{sistema_exponencial}, the Exponential Derivative (derivator function is exponential) of a function $f(x)$ becomes:

\begin{equation} \label{derivadoraExponencial}
\mathfrak{D}\{\}=e^{ax+b}
\end{equation}

\begin{equation} \label{derivadaExponencialParametro_a}
\mathfrak{D}\{a\}\frac{df(x)}{dx}=a^{ins}(x)=\lim_{\Delta \to 0} \frac{ln(f(x+\Delta))-ln(f(x))}{\Delta}
\end{equation}

\begin{equation} \label{derivadaExponencialParametro_b}
\mathfrak{D}\{b\}\frac{df(x)}{dx}=b^{ins}(x)=\lim_{\Delta \to 0} \frac{ln(f(x))(x+\Delta)-ln(f(x+\Delta))x}{\Delta}
\end{equation}

$f(x)$ can be reconstructed from its exponential derivatives as:
\begin{equation} \label{primitivaExponencia}
f(x)=e^{\mathfrak{D}\{a\}\frac{df(x)}{dx}x+\mathfrak{D}\{b\}\frac{df(x)}{dx}}
\end{equation}

The function $e^{-i\omega x}$ (kernel of the Fourier Transform) is discussed in the section \ref{subsec:FourierDerivada}.

\section{Trigonometric Derivators Functions}\label{sec:Trigonometric}
Let $\omega$, $\phi$ $\in \mathbb{R}$, the frequency and phase, respectively, and 
$\omega^{ins}:\mathbb{R} \to \mathbb{C}$, $\phi^{ins}:\mathbb{R} \to \mathbb{C}$ the instantaneous frequency and phase, respectively, the following system can be written:

\begin{equation} \label{sistema_senoidal}
S:\left\{
\begin{array}{ll}
f(x)=\sin(\omega x + \phi)\\
f(x+\Delta)=\sin(\omega (x+\Delta) +\phi)
\end{array}
\right.
\end{equation}

Solving the system and applying the limit of $\Delta \to  0$ in \eqref{sistema_senoidal}, the sinusoidal derivative (derivator function is sinusoidal) of a function $f(x)$ becomes:

\begin{equation} \label{derivadoraSenoidal}
\mathfrak{D}\{\}=\sin({\omega x+\phi)}
\end{equation}

\begin{equation} \label{derivadaSenoidalParametro_w}
\begin{split}
\mathfrak{D}\{\omega\}\frac{df(x)}{dx}=\omega^{ins}(x)=\lim_{\Delta \to 0} \frac{\arcsin(f(x+\Delta))-\arcsin(f(x))}{\Delta} = \\
=\lim_{\Delta \to 0} \frac{\arcsin(f(x+\Delta)\sqrt{1-f(x)^2}-f(x)\sqrt{1-f(x+\Delta)^2})}{\Delta}
\end{split}
\end{equation}

\begin{equation} \label{derivadaSenoidalparametro_phi}
\mathfrak{D}\{\phi\}\frac{df(x)}{dx}=\phi^{ins}(x)=\lim_{\Delta \to 0} \frac{\arcsin(f(x))(x+\Delta)-\arcsin(f(x+\Delta))x}{\Delta}
\end{equation}

The same can be written to cosine and tangent functions:

\begin{equation} \label{derivadoraCosenoidal}
\mathfrak{D}\{\}=\cos({\omega x+\phi)}
\end{equation}

\begin{equation} \label{derivadaCosseinodalParametro_w}
\begin{split}
\mathfrak{D}\{\omega\}\frac{df(x)}{dx}=\omega^{ins}(x)=\lim_{\Delta \to 0} \frac{\arccos(f(x+\Delta))-\arccos(f(x))}{\Delta}=\\
=\lim_{\Delta \to 0} \frac{\arccos(f(x+\Delta)f(x) + \sqrt{(1-f(x+\Delta)^2)(1-f(x)^2)})}{\Delta}
\end{split}
\end{equation}

\begin{equation} \label{derivadaCossenoidalparametro_phi}
\mathfrak{D}\{\phi\}\frac{df(x)}{dx}=\phi^{ins}(x)=\lim_{\Delta \to 0} \frac{\arccos(f(x))(x+\Delta)-\arccos(f(x+\Delta))x}{\Delta}
\end{equation}

\begin{equation} \label{derivadoraTangencial}
\mathfrak{D}\{\}=\tan({\omega x+\phi)}
\end{equation}

\begin{equation} \label{derivadaTangencialParametro_w}
\begin{split}
\mathfrak{D}\{\omega\}\frac{df(x)}{dx}=\omega^{ins}(x)&=\lim_{\Delta \to 0} \frac{\arctan(f(x+\Delta))-\arctan(f(x))}{\Delta}=\\
&=\frac{\arctan(\frac{f(x+\Delta)-f(x)}{1+f(x+\Delta)f(x)})}{\Delta}
\end{split}
\end{equation}

\begin{equation} \label{derivadaTangencialparametro_phi}
\mathfrak{D}\{\phi\}\frac{df(x)}{dx}=\phi^{ins}(x)=\lim_{\Delta \to 0} \frac{\arctan(f(x))(x+\Delta)-\arctan(f(x+\Delta))x}{\Delta}
\end{equation}

If $f(x) \in [-1,1]$, \eqref{derivadaSenoidalParametro_w} to \eqref{derivadaCossenoidalparametro_phi} $\in \mathbb{R}$, otherwise \eqref{derivadaSenoidalParametro_w} to \eqref{derivadaCossenoidalparametro_phi} $\in \mathbb{C}$.

$f(x)$ can be reconstructed from its sinusoidal, cosinusoidal and tangential derivatives, respectively, as:
\begin{equation} \label{primitivaSenoidal}
f(x)=\sin{\biggl(\mathfrak{D}\{\omega\}\frac{df(x)}{dx}x+\mathfrak{D}\{\phi\}\frac{df(x)}{dx}\biggr)}
\end{equation}

\begin{equation} \label{primitivaCossenoidal}
f(x)=\cos{\biggl(\mathfrak{D}\{\omega\}\frac{df(x)}{dx}x+\mathfrak{D}\{\phi\}\frac{df(x)}{dx}\biggr)}
\end{equation}

\begin{equation} \label{primitivaTangencial}
f(x)=\tan{\biggl(\mathfrak{D}\{\omega\}\frac{df(x)}{dx}x+\mathfrak{D}\{\phi\}\frac{df(x)}{dx}\biggr)}
\end{equation}

\section{Instantaneous Frequency and Heisenberg's Uncertainty Principle} \label{sec:Instantaneous}

Let a phase function \cite{Jiang2020} $\Omega(x):\mathbb{R} \to \mathbb{R}$, and a waveform (signal, Wave Function \cite{Born1969}) $\varphi (x,\omega)$ given by:
\begin{equation} \label{waveform}
\varphi (x,\omega)=\sin(\Omega(x))
\end{equation}

If $\Omega(x)$ is known, the determination of the instantaneous frequency $\omega(x)$ presents no difficulty and is determined by:

\begin{equation} \label{derivada_waveform}
\omega(x)=\frac{d\Omega(x)}{dx}
\end{equation}

%its normalized version is:

%\begin{equation} \label{derivada_waveformHZ}
%\omega_{Hz}(x)=\frac{\omega(x)}{2\pi}
%\end{equation}

However, in many real applications, $\Omega(x)$ is not known, but only the waveform $\varphi (x,\omega)$ and then, determining $\omega(x)$ or $x(\omega)$ precisely, from $\varphi (x,\omega)$ is not a possible task, according to Heisenberg's Uncertainty Principle \cite{Heisenberg1927}, \cite{Folland1997}.

Heisenberg's Uncertainty Principle was first proposed for Quantum Mechanics \cite{Ford2005}.
However, it is used to demonstrate that there is a limit to the accuracy with which the pair of canonically conjugate variables \cite{Hjalmars1962} in phase space, ($x,\omega$) or ($x,p$), where $p$ is the momentum, e.g., can be measured simultaneously.\\

\emph{``...Thus, the more precisely the position is determined, the less precisely the momentum is known, and conversely..."} 

\begin{flushright}
\emph{Heisenberg, 1927}  %\cite{Heisenberg1983}   

\end{flushright}
\vspace{.20in}

De Broglie's \cite{LouisVictorPierreRaymond1924} relation establishes the undulatory nature of the particle (matter) by:

\begin{equation} \label{deBroglie}
k=\frac{p}{\hslash}
\end{equation}
where, $k$ is the wavenumber (spatial frequency), $p$ is the momentum and $\hslash$ is the reduced Planck's constant.

Thus, one can understand that determining the momentum $p$ as a function of position $x$ is equivalent to determining the wavenumber $k$ as a function of position $x$ or even the (temporal) frequency $\omega$ as a function of time $x$ ($x$, in this case, is the ``position" in time) - Instantaneous Frequency.

The classical mathematical operation that changes the domain of a function from time to frequency (and vice versa) is the Fourier Transform, but a non-zero function $f(x):\mathbb{R} \to \mathbb{R}$ and its Fourier Transform $F(\omega):\mathbb{R} \to \mathbb{C}$ cannot both be sharply localized \cite{Folland1997}.

\begin{figure}[H]
\centering
\includegraphics[width=0.9\textwidth]{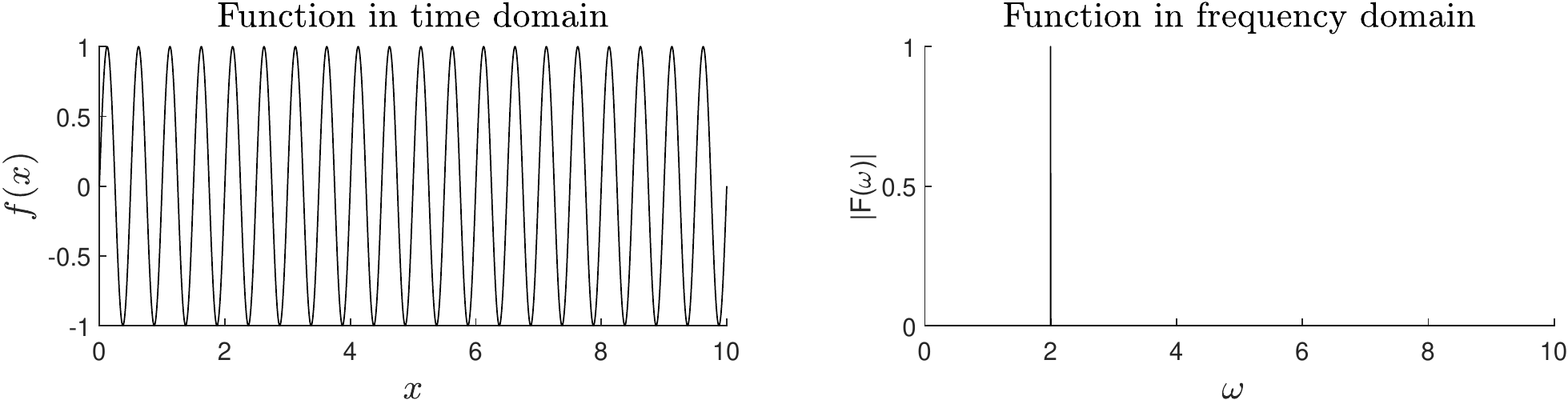}
\caption{Function in time and frequency domain}\label{fig:fourier}
\end{figure}

The figure \ref{fig:fourier} shows a function $f(x):\mathbb{R} \to \mathbb{R}$ (a sinusoidal wave with frequency equal 2 Hz) and $\lvert F(\omega)\rvert:\mathbb{R} \to \mathbb{R}$ in $x$ (time) and $\omega$ (frequency) domain. $F(\omega)$ has no information about $x$ and $f(x)$ has no information about $\omega$.

\begin{figure}[H]
\centering
\includegraphics[width=0.9\textwidth]{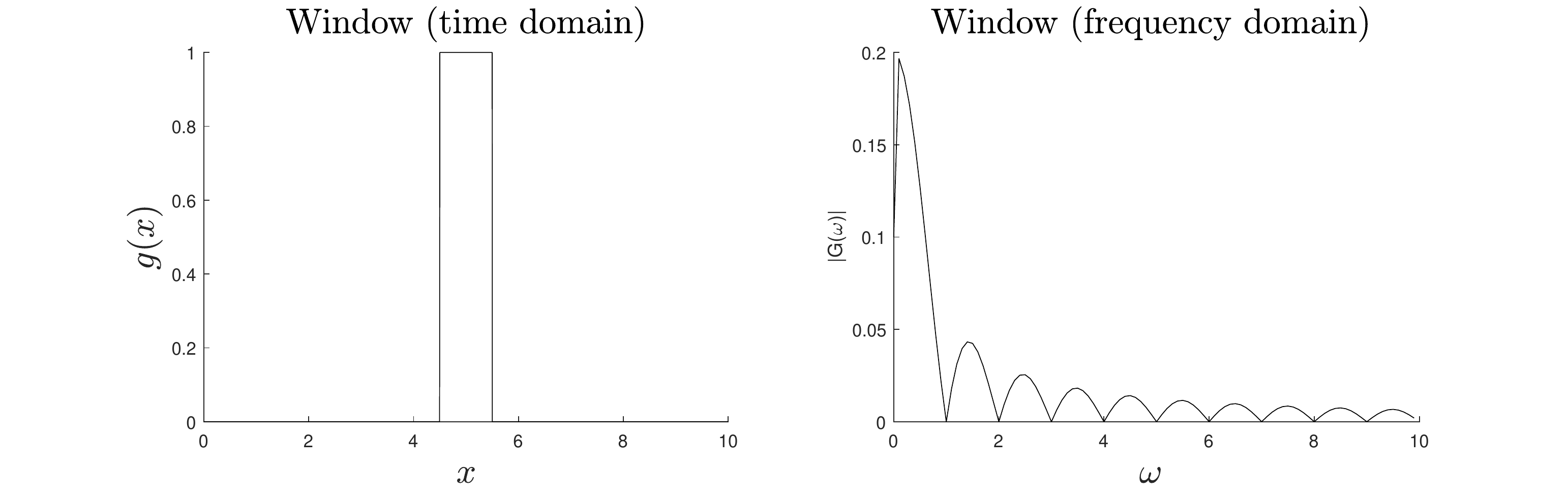}
\includegraphics[width=0.9\textwidth]{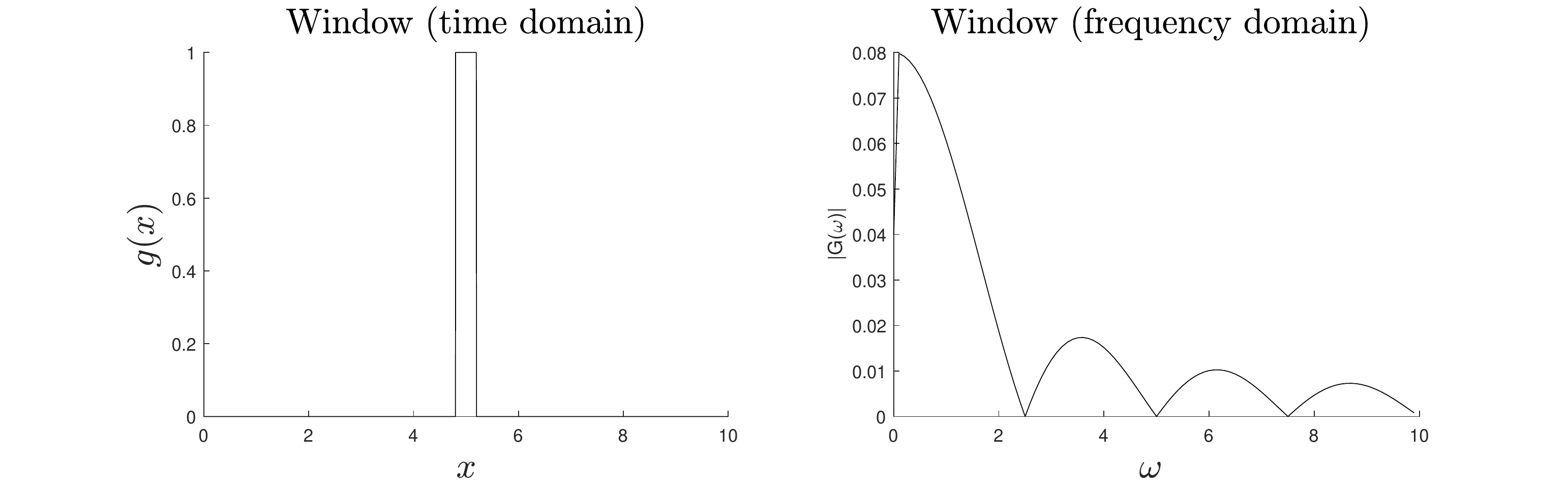}
\caption{Window Functions in time and frequency domain.}
\label{fig:janela}
\end{figure}

A Window Function $g(x):\mathbb{R} \to \mathbb{R}$ that is ``well localized" in the time is used to localize the frequency in time. Figure \ref{fig:janela} shows the wide (above) and narrow (below) window function and its respective Fourier Transform Magnitude  $\lvert G(\omega)\rvert:\mathbb{R} \to \mathbb{R}$ (narrow (above) and wide (below)).

\begin{figure}[H]
\centering
\includegraphics[width=0.9\textwidth]{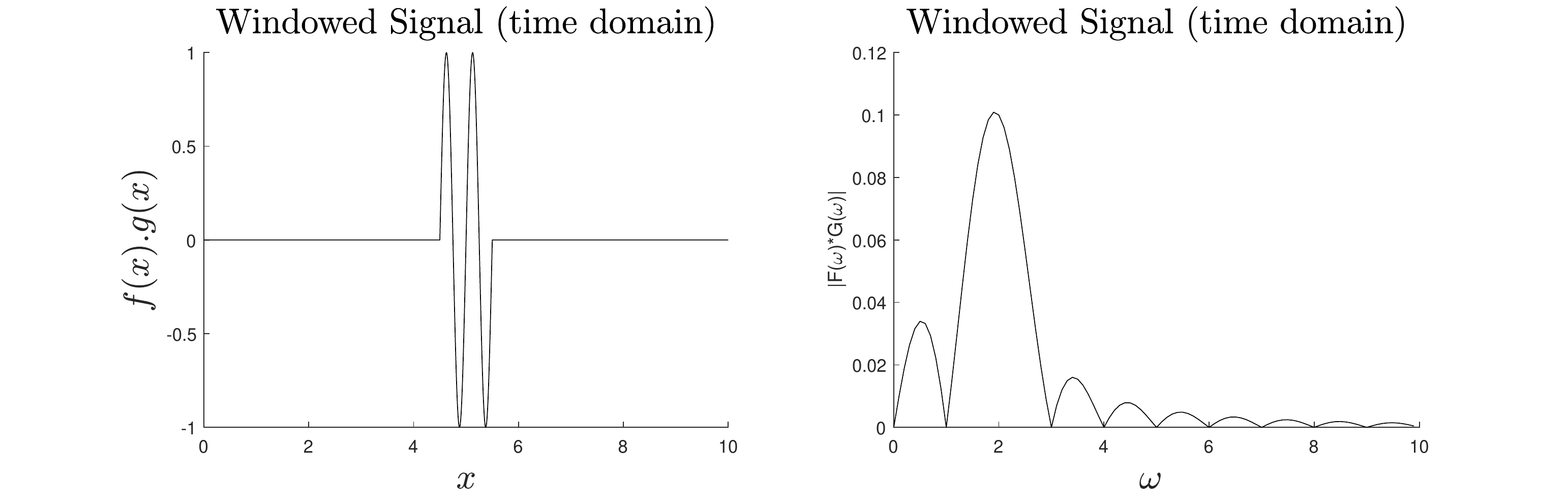}
\includegraphics[width=0.9\textwidth]{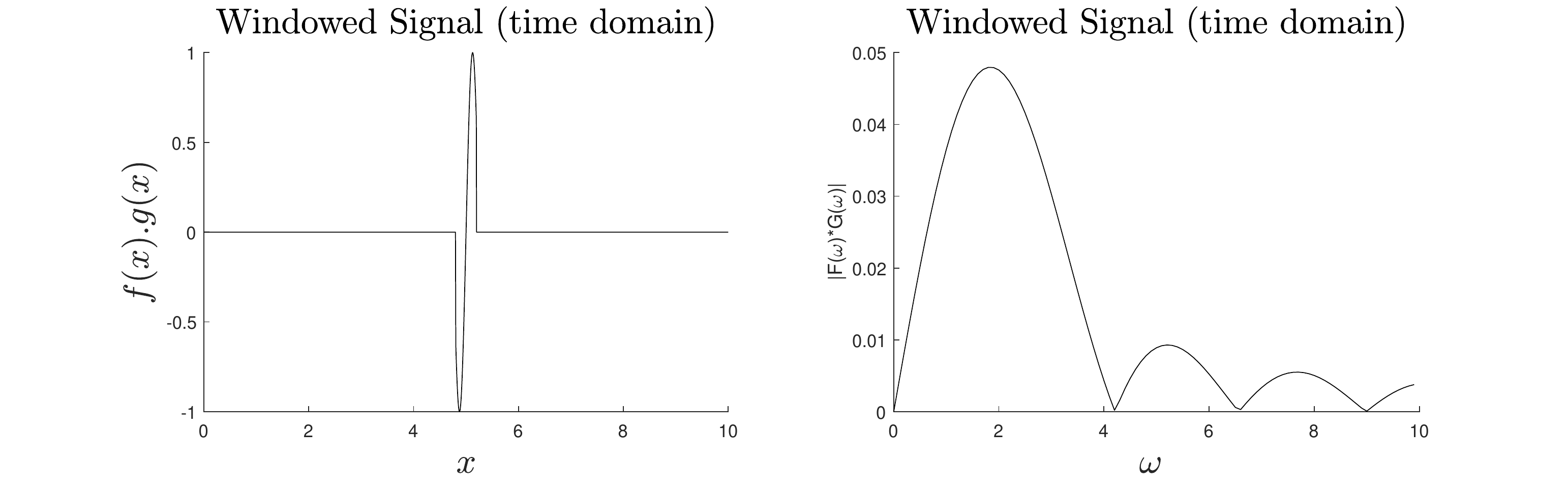}
\caption{Windowed Functions in time and frequency domain.}
\label{fig:sinaljanelado}
\end{figure}

The function $f(x)$ is multiplied by the Window Function $g(x)$. Figure \ref{fig:sinaljanelado} shows the wide (above) and narrow (below) Windowed Function and its respective Fourier Transform Magnitude $\lvert F(\omega)*G(\omega)\rvert$ (narrow (above) and wide (below)).

The Windowed Function in the frequency domain should have only one component (at 2 Hz), but it has components at several frequencies with non-zero amplitudes. This fact is due to the Windowed Function in the frequency domain being the result of the convolution of the function by the Window Function in the frequency domain (in an analogous form, one can use the concept of wave packets in order to locate a wave in space \cite{Rozenman2019}). So it is impossible to identify whether a particular frequency component is due to the function or the Window Function.

To measure the frequency as a function of the time, it was necessary to ``locate" the wave in time using the Window function. However, this fact goes beyond the classical concept in physics of the observer effect \cite{Schlosshauer2005}, \cite{Giacosa2014}, in which to make a measurement, it is necessary to interfere with the measurement (which causes uncertainty).
As everything that exists is a wave (wave nature of matter), Heisenberg's Uncertainty Principle states that uncertainty occurs not only due to the measurement of an experiment (observer effect) but due to the impossibility of locating a wave sharply in the time and frequency (wavenumber, momentum, among others) domain simultaneously.

Analytically, Heisenberg's Uncertainty Principle can be demonstrated considering $\psi(x)$ and $\Psi(p)$ wave functions and Fourier Transform \footnote{$\psi(x)$ and $\Psi(p)$ are functions in two corresponding orthonormal bases in Hilbert space and, therefore, are Fourier Transform of each other and $x$ and $p$ are conjugate variables.} of each other for position $x$ and momentum $p$ , respectively. 

Born's rule \cite{Born1926} states that $\lvert\psi(x)\rvert^2$ and $\lvert\Psi(p)\rvert^2$ are probability density functions and then the variances of position $\sigma_x^2$ and momentum $\sigma_p^2$ are:

\begin{equation} \label{varianciaX}
\sigma_x^2=\int_{-\infty}^{\infty}{x^2}\lvert\psi(x)\rvert^2\,dx\ = 
\int_{-\infty}^{\infty}{(x\psi(x))^*x\psi(x)}\,dx\
\end{equation} 

\begin{equation} \label{varianciaP}
\sigma_p^2=\int_{-\infty}^{\infty}{p^2}\lvert\Psi(p)\rvert^2\,dp\ =
\int_{-\infty}^{\infty}{(p\Psi(p))^*p\Psi(p)}\,dp\
\end{equation} 

Let $f(x)=x\psi(x)$:
\begin{equation} \label{varianciaXInner}
\sigma_x^2=\int_{-\infty}^{\infty}{f^*(x)f(x)}\,dx\ =
\int_{-\infty}^{\infty}{\lvert f(x)\rvert^2}\,dx\ =
\langle f \lvert f \rangle
\end{equation} 

Let $\mathcal{F} \left\{.\right\}$ the Fourier transform, $-i{\hslash}\frac{d\psi(x)}{dx}$ the momentum operator in position space,  $G(p)=p\Psi(p)$, $g(x)=\mathcal{F} \left\{G(p)\right\}$ and applying the Parseval's theorem \cite{Stone2021}:
\begin{equation} \label{varianciaPInner}
\sigma_p^2=\int_{-\infty}^{\infty}{G^*(p)G(p)}\,dp\ =
\int_{-\infty}^{\infty}{\lvert G(p)\rvert^2}\,dp\ =
\int_{-\infty}^{\infty}{\lvert g(x)\rvert^2}\,dx\ =
\langle g \lvert g \rangle
\end{equation} 

Using the Cauchy–Schwarz inequality \cite{Zhou2021}:
\begin{equation} \label{Cauchy–Schwarz}
\langle f \lvert f \rangle \langle g \lvert g \rangle \geq \lvert \langle f \lvert g \rangle \rvert^2 
\end{equation} 

\begin{equation} \label{Inner(f,g)}
\lvert \langle f \lvert g \rangle \rvert^2  \geq \operatorname{Im}(\lvert \langle f \lvert g \rangle \rvert^2 ) = \biggl(\frac{\langle f \lvert g \rangle - \langle g \lvert f \rangle}{2i}\biggr)^2
\end{equation} 

\begin{equation} \label{carai}
\begin{split}
\langle f \lvert g \rangle &- \langle g \lvert f \rangle = \\
=\int_{-\infty}^{\infty}{x\psi^*(x)\biggl(-i{\hslash}\frac{d\psi(x)}{dx}\biggr)}\,dx\ &-
\int_{-\infty}^{\infty}{\biggl(-i{\hslash}\frac{d\psi^*(x)}{dx}\biggr)x\psi(x)}\,dx\ =
i{\hslash}
\end{split}
\end{equation} 

Applying \ref{carai} in \ref{Inner(f,g)}:

\begin{equation} \label{Inner(f,g)Quadrado}
\lvert \langle f \lvert g \rangle \rvert^2 = \biggl(\frac{i{\hslash}}{2i}\biggr)^2 =
\frac{{\hslash^2}}{4}
\end{equation}

Applying \ref{varianciaXInner}, \ref{varianciaPInner} and \ref{Inner(f,g)Quadrado} in \ref{Cauchy–Schwarz}:

\begin{equation} \label{PrincipioIncertezaSigma}
\sigma_x^2\sigma_p^2 \geq \frac{{\hslash^2}}{4} = \sigma_x\sigma_p \geq \frac{{\hslash}}{2}
\end{equation} 

The demonstration of the uncertainty principle is strictly mathematical. Any pair of variables conjugated will produce the same results as this demonstration.

Following Kennard's consideration \cite{Choe2020}, $\Delta x = \sigma_x$ the uncertainty in position $x$ (proportional to the width of the Window Function in time or space domain), $\Delta p = \sigma_p$ the uncertainty in momentum, $h$ the Plank's constant, Heisenberg’s Uncertainty Principle is normally presented as:

\begin{equation} \label{PrincipioIncertezaFINAL}
\Delta x \Delta p \geq \frac{h}{4\pi}
\end{equation} 

In the frequency domain, let $\Delta k$ be the uncertainty in wavenumber (spatial frequency) or $\Delta\omega$ the uncertainty in (temporal) frequency ($\Delta k$ or $\Delta\omega$ are proportional to the width of the Window Function in frequency domain). Through de Broglie´s relation \ref{deBroglie}, \ref{PrincipioIncertezaFINAL} can be written as: 
\begin{equation} \label{PrincipioIncerteza}
\Delta x \Delta k \geq \frac{1}{2\pi} \:\:\:\:  \textrm{or}\:\:\:\:  \Delta x \Delta \omega \geq \frac{1}{2\pi}  
\end{equation}

Another way to understand Heisenberg’s Uncertainty Principle (and perhaps the simplest) is through the Fourier Transform of the Gaussian function.
Let $f(x)$ be a Gaussian function in the space (time) domain $x$, $F(\omega)$ is its Fourier transform in the frequency domain $\omega$ and is also a Gaussian function. Then, the standard deviation $\sigma$ can be understood as a measure of precision, and this occurs inversely in $f(x)$ and $F(\omega)$. Thus, if the uncertainty is small in one domain, it is large in the other domain.

\begin{equation} \label{fourierGaussiana}
f(x)=e^{-\sigma x^2} \Leftrightarrow F(\omega)=\frac{1}{\sqrt{2\sigma}}e^{-\frac{\omega^2}{4\sigma}}
\end{equation} 

A time-frequency representation\footnote{Time-frequency is a representation with a two-dimensional domain $(x,\omega)$, and is used to represent any pair of canonically conjugate coordinates, such as time-frequency, position-wavenumber, position-momentum, among others.} is used when it is necessary to ``localize" $\omega$ in $x$ (instantaneous frequency) and vice-versa. This representation is also known as a spectrogram, generally obtained through the Short Time Fourier Transform \cite{Groechenig2001} (or other transforms, such as  Wavelet Transform \cite{Strang1996}, Wigner-Ville distribution function \cite{Hlawatsch1997}, etc).

A classic example of a signal whose frequency varies with time (non-stationary signal) is the Chirp Signal \cite{Gretinger2014}.

% In the Short Time Fourier Transform, the Fourier Transform is applied to the signal multiplied by a window function in order to also locate, in time, the representation of the signal in the frequency domain. REFERENCIAS

\subsection{Chirp Derivative}\label{subsec:ChirpDerivada}

A Chirp Signal can be defined with the following waveform (signal):
\begin{equation} \label{chirp_seno}
f(x)=\sin(2\pi\int{\omega(x)}dx+\phi),
\end{equation}
where, $\omega(x)$ is the instantaneous frequency function and $\Omega(x)=2\pi\int{\omega(x)}dx+\phi$ is the phase function.\\

The following example is suitable for showing the determination of the instantaneous frequency. Considering a Chirp Signal with instantaneous frequency function given by:

\begin{equation} \label{frequencia_instantanea_quadratica_exemplo}
\omega(x)=x^2+2x+1
\end{equation}

For $\phi=0$, the waveform is;

\begin{equation} \label{chirp_quadratico_exemplo}
f(x)=\sin \biggl(2\pi\biggl(\int_{}^{}({x^2+2x+1}) \,dx \biggr) + 0 \biggr) = \sin\biggl(2\pi\biggl(\frac{x^3}{3}+x^2+x\biggr)\biggr)
\end{equation}

\begin{figure}[H]
\centering
\includegraphics[width=0.9\textwidth]{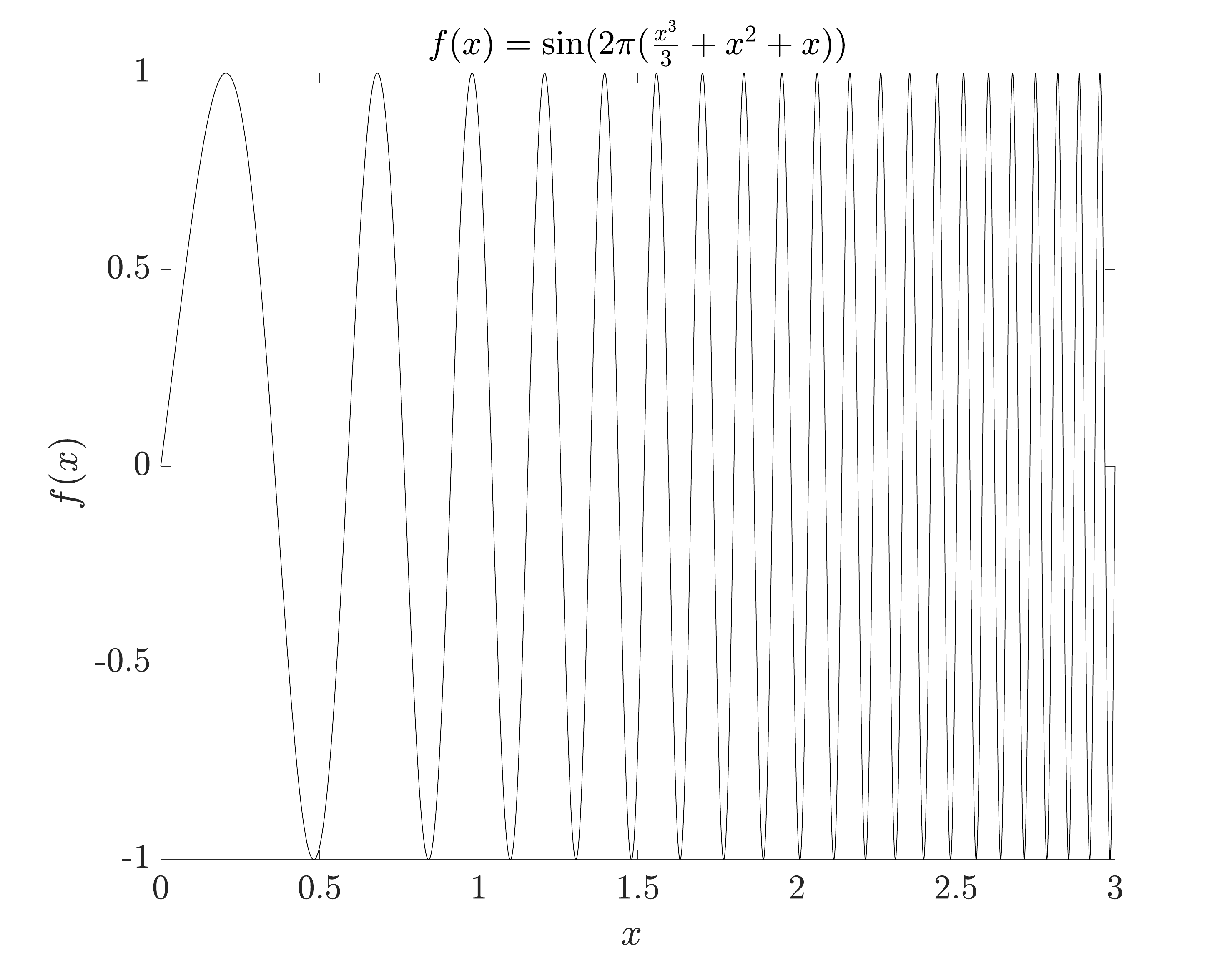}
\caption{Quadratic Chirp Signal - $\omega(x)=x^2+2x+1$;  $f(x)=\sin\biggl(2\pi\biggl(\frac{x^3}{3}+x^2+x\biggr)\biggr)$}\label{fig4}
\end{figure}

\begin{figure}[H]
\centering
\includegraphics[width=0.9\textwidth]{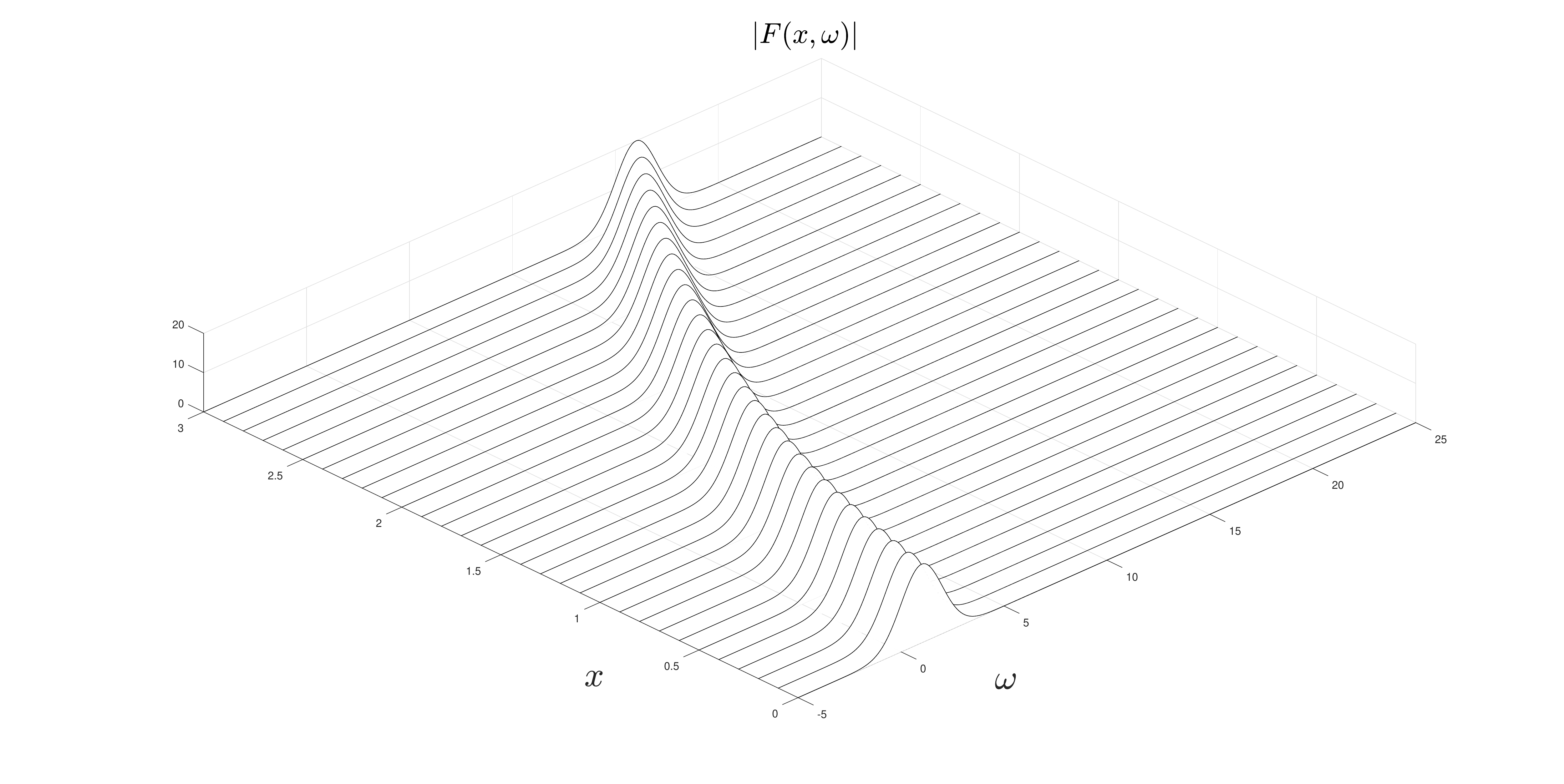}
\caption{Short Time Fourier Transform - $\lvert F(x,\omega) \rvert$}\label{figSTFT}
\end{figure}

Figure \ref{figSTFT} shows the spectrogram for the signal sampled at 100 samples/unit of $x$ (100Hz if $x$ is given in seconds) obtained through Short Time Fourier Transform with Gaussian window function (Gabor Transform \cite{Groechenig2001}) and standard deviation equal to 1. The $\lvert F(x,\omega) \rvert$  values (z-axis) are proportional to the energy in the signal at $(x,\omega)$. For each $x$ (or $\omega$) value there is a range of $\omega$ (or $x$) values whose function $\lvert F(x,\omega) \rvert$ is non-zero. These intervals at $(x,\omega)$ domain represent the uncertainty in determining the instantaneous frequency in this signal representation.

%\begin{figure}[H]
% \centering
% \includegraphics[width=0.9\textwidth]{stft.eps}
% \caption{STFT - $\Omega_q(x)=x^2+2x+1$,  $f(x)=\sin\biggl(2\pi\biggl(\frac{x^3}{3}+x^2+x\biggr)+\phi\biggr)$}\label{fig5}
% \end{figure}

With:
\begin{equation} \label{frquencia linear}
\omega(x)=\omega_1 x+\omega_0
\end{equation}
\eqref{chirp_linear} is a Linear Chirp (or Quadratic Phase Signal), with initial frequency $\omega_0$ (at $x=0$) and rate chirp $\omega_1$:
\begin{equation} \label{chirp_linear}
f(x)=\sin \biggl(2\pi\biggl(\int_{}^{}({\omega_1 x+\omega_0}) \,dx \biggr) +\phi \biggr) = \sin\biggl(2\pi\biggl(\frac{1}{2}\omega_1 x^2+\omega_0 x\biggr)+\phi\biggr)
\end{equation}

The resolution of the system,

\begin{equation} \label{sistema_chirp}
S:\left\{
\begin{array}{ll}
f(x)=\sin\biggl(2\pi\biggl(\frac{1}{2}\omega_1 x^2+\omega_0 x\biggr)+\phi\biggr)\\
f(x+\Delta)=\sin\biggl(2\pi\biggl(\frac{1}{2}\omega_1 (x+\Delta)^2+\omega_0 (x+\Delta)\biggr)+\phi\biggr)\\
f(x+2\Delta)=\sin\biggl(2\pi\biggl(\frac{1}{2}\omega_1 (x+2\Delta)^2+\omega_0 (x+2\Delta)\biggr)+\phi\biggr)\\
\end{array}
\right.
\end{equation}

applying limit to $\Delta\to 0$, is:
%\begin{eqnarray} \label{w_solucao_sistema_chirp}
%\mathfrak{D}\{\omega_1\}\frac{df(x)}{dx} =\omega_1(x)= %\;\;\;\;\;\;\;\;\;\;\;\;\;\;\;\;\;\;\;\;\;\;\;\;\;\;\;\;\;\;\;\;\;\;\;\;\;\;\;\;\;\;\;\;\;\;%\;\;\;\;\;\;\;\;\;\;\;\; & & \nonumber \\
%\lim_{\Delta \to 0} \frac{\arcsin(f(x))-\arcsin(2f(x+\Delta))+\arcsin(f(x+2\Delta))}{2\pi %\Delta^2}
%\end{eqnarray}

\begin{equation} \label{w_solucao_sistema_chirp} 
\mathfrak{D}\{\omega_1\}\frac{df(x)}{dx}=\omega_1^{ins}(x)=       
\lim_{\Delta \to 0} \frac{K_0-K_1+K_2}{2\pi \Delta^2}\\
\end{equation}

where, 
\begin{align}
K_0 &= \arcsin(f(x)) \nonumber \\
K_1 &= \arcsin(2f(x+\Delta))\nonumber \\
K_2 &= \arcsin(f(x+2\Delta))\nonumber 
\end{align}

%\begin{eqnarray} \label{w0_solucao_sistema_chirp}
%\mathfrak{D}\{\omega_0\}\frac{df(x)}{dx}=\omega_0(x)= 
%\lim_{\Delta \to 0} -\frac{1}{2} \frac{(2x+3\Delta)\arcsin(f(x))-%(4x+4\Delta)\arcsin(f(x+\Delta))+(2x+\Delta)\arcsin(f(x+2\Delta))}{2\pi \Delta^2}
%\end{eqnarray}

\begin{equation} \label{w0_solucao_sistema_chirp} 
\mathfrak{D}\{\omega_0\}\frac{df(x)}{dx}=\omega_0^{ins}(x)=       
\lim_{\Delta \to 0} -\frac{1}{2} \frac{K_0-K_1+K_2}{2\pi \Delta^2}\\
\end{equation}

where, 
\begin{align}
K_0 &= (2x+3\Delta)\arcsin(f(x)) \nonumber \\
K_1 &= (4x+4\Delta)\arcsin(f(x+\Delta))\nonumber \\
K_2 &= (2x+\Delta)\arcsin(f(x+2\Delta))\nonumber 
\end{align}

%\begin{eqnarray} \label{fase_solucao_sistema_chirp1}
%\mathfrak{D}\{\phi\}\frac{df(x)}{dx}=\phi(x) =\;\;\; %\;\;\;\;\;\;\;\;\;\;\;\;\;\;\;\;\;\;\;\;\;\;\;\;\;\;\;\;\;\;\;\;\;\;\;\;\;\;\;\;\;\;\;\;\;\;%\;\;\;\;\;\;\;\;\;\;\;\;\;\;\;\;\;\;\;\;\;\;\;\;\;\;\;\;\;\;\;\;\;\ & & \nonumber \\
 %\lim_{\Delta \to 0} \frac{1}{2} \frac{(2\Delta^2+x^2+3x\Delta)\arcsin(f(x))+(-2x^2-%4x\Delta)\arcsin(f(x+\Delta))+(x^2+x\Delta)\arcsin(f(x+2\Delta))}{\Delta^2}\;\;\;\;\;\; & &
%\end{eqnarray}

\begin{equation} \label{fase_solucao_sistema_chirp} 
\mathfrak{D}\{\phi\}\frac{df(x)}{dx}=\phi(x)=       
\lim_{\Delta \to 0} \frac{1}{2} \frac{K_0+K_1+K_2}{\Delta^2}
\end{equation}
where, 
\begin{align}
K_0 &= (2\Delta^2+x^2+3x\Delta)\arcsin(f(x))\nonumber \\
K_1 &= (-2x^2-4x\Delta)\arcsin(f(x+\Delta))\nonumber \\
K_2 &= (x^2+x\Delta)\arcsin(f(x+2\Delta))\nonumber 
\end{align}

According to \eqref{frquencia linear}, the instantaneous frequency function is obtained as:
\begin{equation} \label{frquencia linear reconstruida}
\omega(x)=\mathfrak{D}\{\omega1\}\frac{df(x)}{dx}x+\mathfrak{D}\{\omega_0\}\frac{df(x)}{dx}
\end{equation}

$\mathfrak{D}\{\omega_1\}$, $\mathfrak{D}\{\omega_0\}$ ($\mathfrak{D}\{\phi\}$ is not needed in this example) are:
\begin{equation} \label{w1_solucao_chirp_quadratico_exemplo}
\mathfrak{D}\{\omega_1\}\frac{df(x)}{dx} =\frac{2\, {\cos\!\left(\frac{2\, \pi\, x\, \left(x^2 + 3\, x + 3\right)}{3}\right)}^3\, \left(x + 1\right)}{{\left({\cos\!\left(\frac{2\, \pi\, x\, \left(x^2 + 3\, x + 3\right)}{3}\right)}^2\right)}^{\frac{3}{2}}}
\end{equation}

\begin{equation} \label{w0_solucao_chirp_quadratico_exemplo}
\mathfrak{D}\{\omega_0\}\frac{df(x)}{dx} =-\frac{{\cos\!\left(\frac{2\, \pi\, x\, \left(x^2 + 3\, x + 3\right)}{3}\right)}^3\, \left(x^2 - 1\right)}{{\left({\cos\!\left(\frac{2\, \pi\, x\, \left(x^2 + 3\, x + 3\right)}{3}\right)}^2\right)}^{\frac{3}{2}}}
\end{equation}

% \begin{equation} \label{fase_solucao_chirp_quadratico_exemplo}
% \mathfrak{D}\{\phi\}\frac{df(x)}{dx} =-\frac{4\, \pi\, x\, {\cos\!\left(\frac{2\, \pi\, x\, \left(x^2 + 3\, x + 3\right)}{3}\right)}^3 - 2\, \arcsin\!\left(\sin\!\left(\frac{2\, \pi\, x\, \left(x^2 + 3\, x + 3\right)}{3}\right)\right)\, {\left({\cos\!\left(\frac{2\, \pi\, x\, \left(x^2 + 3\, x + 3\right)}{3}\right)}^2\right)}^{\frac{3}{2}} + 4\, \pi\, x^2\, {\cos\!\left(\frac{2\, \pi\, x\, \left(x^2 + 3\, x + 3\right)}{3}\right)}^3}{4\, \pi\, {\left({\cos\!\left(\frac{2\, \pi\, x\, \left(x^2 + 3\, x + 3\right)}{3}\right)}^2\right)}^{\frac{3}{2}}}
% \end{equation}
%ESSA DEFERIVADA DA FASE NÃO PRECISA AQUI, NÃO USA PARA NADA NESTE EXEMPLO

The \eqref{w1_solucao_chirp_quadratico_exemplo} and \eqref{w0_solucao_chirp_quadratico_exemplo} have positive and negative values, which are therefore associated with positive and negative frequency values. In absolute values, \eqref{w1_solucao_chirp_quadratico_exemplo} and \eqref{w0_solucao_chirp_quadratico_exemplo} are, respectively:

\begin{equation} \label{Modulo_w1_solucao_chirp_quadratico_exemplo}
\left\lvert{\mathfrak{D}\{\omega_1\}\frac{df(x)}{dx}}\right\rvert = \lvert 2x+2 \rvert
\end{equation}

\begin{equation} \label{Modulo_w0_solucao_chirp_quadratico_exemplo}
\left\lvert{\mathfrak{D}\{\omega_0\}\frac{df(x)}{dx}}\right\rvert = \lvert -x^2+1 \rvert
\end{equation}

However, negative frequency values can be neglected, and therefore the modulus functions at \eqref{Modulo_w1_solucao_chirp_quadratico_exemplo} and \eqref{Modulo_w0_solucao_chirp_quadratico_exemplo} can be removed without loss of generality. According to \eqref{frquencia linear reconstruida}, the frequency function $\omega_{qr}(x)$ can be reconstructed by:
\begin{equation} \label{frquencia linear reconstruida exemplo}
\omega_{qr}(x)=(2x+2)x-x^2+1 = x^2+2x+1
\end{equation}

\begin{figure}[H]
\centering
\includegraphics[width=0.9\textwidth]{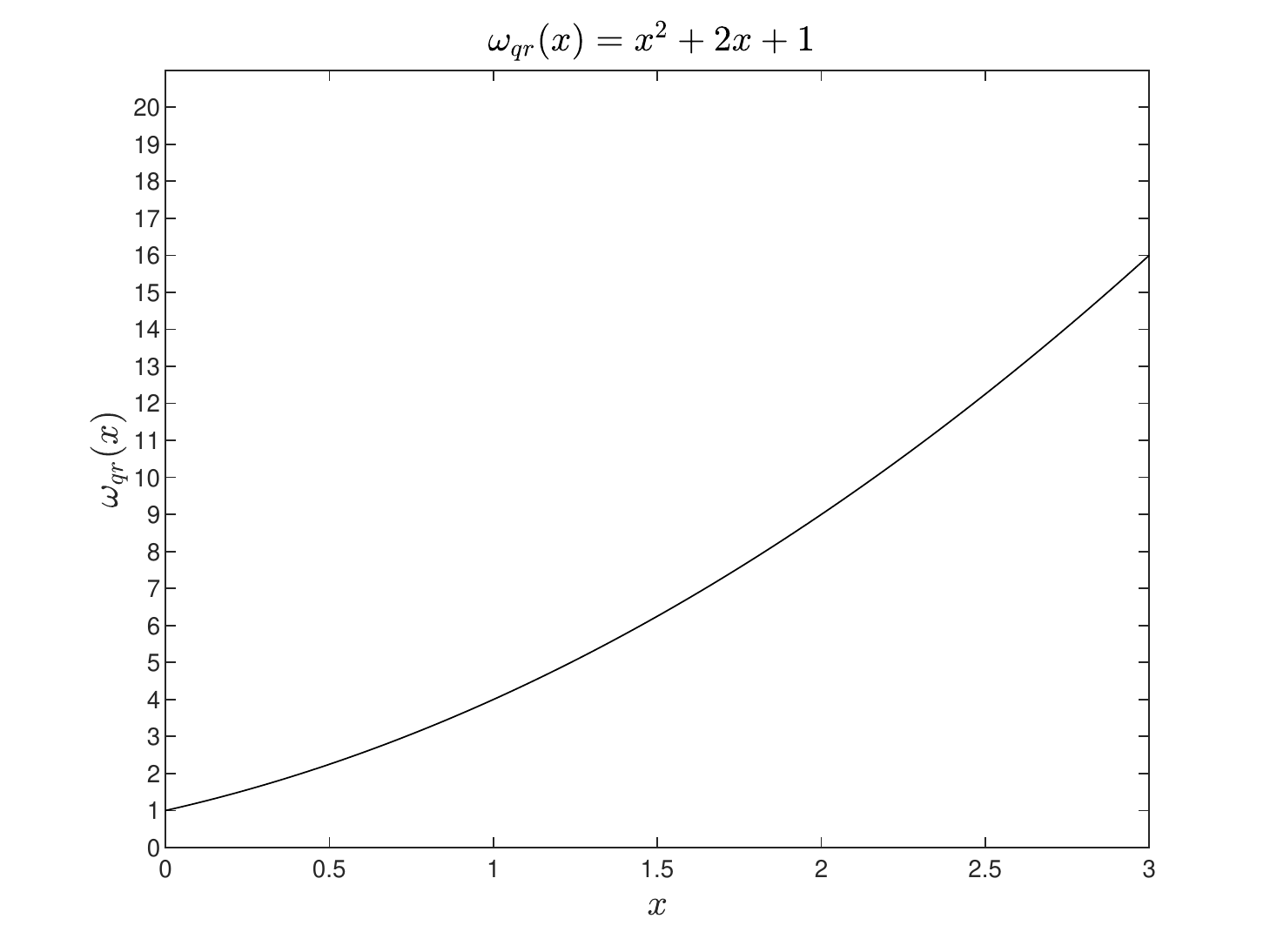}
\caption{$\omega_{qr}(x)=x^2+2x+1$, obtained from $f(x)=\sin\biggl(2\pi\biggl(\frac{x^3}{3}+x^2+x\biggr)+\phi\biggr)$}\label{figWQR}
\end{figure}

It is important to note that \eqref{frquencia linear reconstruida exemplo} is obtained from \eqref{chirp_quadratico_exemplo} and not from the \eqref{derivada_waveform}, i.e., the exact instantaneous frequency is obtained from waveform (wave function or signal) ($f(x)$) and not from phase function ($\Omega(x))$ and \textbf{there is no uncertainty.} 

\begin{remark}
Heisenberg’s Uncertainty Principle was not respected.    
\end{remark}

% and, acoording to \eqref{chirp_linear}, $f(x)$ as:
%  \begin{eqnarray} \label{primitivaChirp}
%  f(x)&=&\sin \Biggl\{ 2\pi\left(\frac{1}{2}\mathfrak{D}\{\omega1\}\frac{df(x)}{dx}x^2+\mathfrak{D}\{\omega_0\}\frac{df(x)}{dx}x \right)  +\mathfrak{D}\{\phi\}\frac{df(x)}{dx}  \Biggl\} 
%  \end{eqnarray}
% ESSE COMENTÁRIO NÃO ESTÁ DANDO CERTO. PRECISA PENSAR MELHOR.
% MAS PODE DEIXAR SEM ISSO

From the instantaneous frequency function $\omega(x)$, the amplitude spectrum $F(\omega)$ can be calculated as:
\begin{equation} \label{histograma_frequencia_instantanea}
F(\omega)=\frac{1}{\omega} \int_{-\infty}^{\infty} \omega(x)h(x)dx,
\end{equation}
with,
\begin{equation} \label{janela_histograma_frequencia_instantanea}
h(x) = \begin{cases} 1;\;\textrm{ if } \omega(x) = \omega;\\ 0;\;\textrm{ otherwise }\\ \end{cases}
\end{equation}

\subsection{Fourier Derivative}\label{subsec:FourierDerivada}

An important exponential function is $e^{-i\omega x}$, with $i=\sqrt{-1}$ and $x,\omega, \in \mathbb{R}$, which is the kernel of the Fourier Transforms \cite{Stone2021}. Adding a scaling factor (as in \eqref{exponencial1}) to this kernel, and proceeding analogously to \eqref{sistema_exponencial}, the \textbf{Fourier Derivative} of function $f(x)$ is:

\begin{equation} \label{derivadoraFourier}
\mathfrak{D}\{\}=e^{-i\omega x + b}
\end{equation}

\begin{equation} \label{derivadaFourierParametro_w}
\mathfrak{D}\{\omega\}\frac{df(x)}{dx}=\omega(x)=\lim_{\Delta \to 0} \frac{ln(f(x+\Delta))-ln(f(x))}{-i\Delta}
\end{equation}

\begin{equation} \label{derivadaFourierParametro_b}
\mathfrak{D}\{b\}\frac{df(x)}{dx}=b(x)=\lim_{\Delta \to 0} \frac{ln(f(x))(x+\Delta)-ln(f(x+\Delta))x}{\Delta}
\end{equation}

$f(x)$ can be reconstructed from its Fourier derivatives as:
\begin{equation} \label{primitivaFouier}
f(x)=e^{-i\mathfrak{D}\{\omega\}\frac{df(x)}{dx}x+\mathfrak{D}\{b\}\frac{df(x)}{dx}}
\end{equation}

The parameter $\omega$ is the frequency in the kernel of the Fourier Transform, and therefore \eqref{derivadaFourierParametro_w} is the instantaneous frequency ($w(x) \in \mathbb{C}$). 

Let $A \in \mathbb{R}$, a wave function of the type $\psi(x,\omega)=Ae^{-i \omega x}$. Considering, as example, $A=2$, $\omega(x) = x^3+2x$ (the frequency $\omega$ varies with $x$, i.e. $\omega(x)$) and $\Omega(x)$ the phase function (as in Chirp Signal \eqref{chirp_seno}), the wave function $\psi(x,\omega)$ becomes:

\begin{equation} \label{wavefunctionExample}
\psi(x) = 2e^{-i \Omega(x)} = 2e^{-i\left(\frac{x^4}{4}+x^2\right)}
\end{equation}
where, 
\begin{equation} \label{wavefunctionExample_phasefunction}
\Omega(x) = \int_{}^{}(x^3+2x) dx = \left(\frac{x^4}{4}+x^2\right)
\end{equation}

Applying \ref{derivadaFourierParametro_w} and \ref{derivadaFourierParametro_b} in \ref{wavefunctionExample}:

\begin{equation} \label{derivadaFourierParametro_w_WaveFunction}
\mathfrak{D}\{\omega\}\frac{d\psi(x)}{dx}=\omega(x)=x^3+2x
\end{equation}

\begin{equation} \label{derivadaFourierParametro_b_WaveFucntion}
\mathfrak{D}\{b\}\frac{d\psi(x)}{dx}=b(x)=i\left(\frac{3x^4}{4}+x^2\right)+\ln{2}
\end{equation}

Applying \ref{primitivaFouier} in \ref{derivadaFourierParametro_w_WaveFunction} and \ref{derivadaFourierParametro_b_WaveFucntion}, the wave function $\psi(x)$ \eqref{wavefunctionExample}, can be reconstructed as:

\begin{equation} \label{reconstruida_wave_function}
\psi(x)=e^{-i(x^3+2x)x+i\left(\frac{3x^4}{4}+x^2\right)+\ln{2}}=2e^{-i\left(\frac{x^4}{4}+x^2\right)}
\end{equation}

\textbf{There is no uncertainty.} 

\begin{remark}
Heisenberg’s Uncertainty Principle was not respected.    
\end{remark}

\section{Conclusion}\label{sec:Conclusion}

In a simplified and summarized manner, Differential Calculus is based on applying a limit tending to zero for Newton’s Difference Quotient applied under any function $f(x)$. 

This operation determines another function (the derivative) whose values represent the instantaneous angular coefficients of the tangent lines to the function $f(x)$.

This paper showed that the Differential and Integral Calculus could be applied to other parameters of other functions called derivator and integrator functions.

All the theories presented can be applied to two or more dimensions (partial derivatives and multiple integrals), in addition to well-established operations in classical differential and integral calculus such as the chain rule, product and division derivatives and integrals, differential and integral equations, and others, and this is suggested as future work.

Some examples were presented, with emphasis on the determination of the instantaneous frequency.
Although Heisenberg's Uncertainty Principle is formalized as a property of waves, this paper has shown that uncertainty occurs due to the methodology employed for determining the instantaneous frequency in a function (wave function or signal). 

Heisenberg’s Uncertainty Principle is based on the use of integral transforms (such as
Fourier Transform and similar wave packets), for a function in the time (or space)
domain to obtain its representation in the frequency domain and vice versa.

An integral transform is obviously based on the calculation of integrals. Hence, the integral is suitable for measuring general quantities associated with the whole function domain, such as an area, expected value, norm, autocorrelation, and even frequency distribution (spectral density), but not instantaneous quantities.

Integral transforms (or wave packets) will produce uncertainty in the phase space of canonically conjugate variables.

Nevertheless, why use a mathematical operation based on integral to try to determine instantaneous quantities?

In turn, the derivative is suitable for measuring instantaneous
quantities in a function. This paper presented a form to obtain the instantaneous frequency of a function given in the time (or space) domain using
derivatives (and not integrals).

The Fourier, Trigonometric, and Chirp Derivatives are examples of different forms to obtain the instantaneous frequency sharply.

%Heisenberg’s Uncertainty Principle no longer applies in the form it was conceived, being reduced to the observer effect when there is some uncertainty in some experiment/measurement.

%In brief, the uncertainty is due solely to the methodology and/or the instrument used to make the measurement and not to an inherent property of matter or waves.

%\begin{flushright}
%``God does not play dice" (Albert Einstein, 1926).
%\end{flushright}

%\emph{``...God does not play dice..."} 

%\begin{flushright}
%\emph{Albert Einstein, 1926}  

%\end{flushright}

\pagebreak

% \begin{equation}
% \mbox{\Huge $\supset \!\!\!\!\!\! \vert  \mbox{\small \raisebox{4pt}{a}} $}\; : (x,y) = \frac{a}{b}
% \end{equation}

% \begin{equation}
% \mbox{\Huge $\supset \!\!\!\!\!\! \vert  \mbox{\small \raisebox{5pt}{\underline{a}}} $}\; : (x,y) = \frac{a}{b}
% \end{equation}

% \begin{equation}
% \mbox{\Huge $\supset \!\!\!\!\!\! \vert  \mbox{\small \raisebox{5pt}{\underline{{\it a}}}} $}\; : (x,y) = \frac{a}{b}
% \end{equation}

%\derivadora{p} \;$: ax+b$

%{\Huge$\supset\!\!\!\!\!\vert$}\!\!\!$^a$

%$\supset\!\!\!\!\!\vert$ \!\!\!$^a$

%\bibliography{teste.bib}% common bib file   

\end{document}